\documentclass[11pt]{article}
\usepackage{amsmath,amsfonts}
\usepackage{amsthm,amssymb}
\usepackage{epsfig}
\usepackage{leftidx}
\usepackage{graphicx}
\usepackage{enumitem}
\usepackage{color, epstopdf}
\usepackage{cite}
\usepackage{graphicx,wrapfig,subfigure}

\newtheorem{theorem}{Theorem}[section]
\newtheorem{lemma}{Lemma}[section]

\newtheorem{example}{Example}
\newtheorem{remark}{Remark}

\newcommand{\dfd}[1]{\partial^{#1}_\tau}

\newcommand{\diff}{\triangledown_\tau}
\newcommand{\piff}{\partial_\tau}
\newcommand{\defeq}{:=}
\newcommand{\zd}{\,\mathrm{d}}
\newcommand{\abs}[1]{\left|#1\right|}

\newcommand{\bra}[1]{\left(#1\right)}
\newcommand{\brab}[1]{\big(#1\big)}
\newcommand{\braB}[1]{\Big(#1\Big)}
\newcommand{\brat}[1]{(#1)}
\newcommand{\kbra}[1]{\left[#1\right]}

\newcommand{\myinner}[1]{\left\langle#1\right\rangle}
\newcommand{\myinnerb}[1]{\big\langle#1\big\rangle}

\newcommand{\mynorm}[1]{\left\|#1\right\|}
\newcommand{\mynormb}[1]{\big\|#1\big\|}

\def\lan#1{\textcolor{blue}{#1}}

\begin{document}
\title{Asymptotically compatible energy and dissipation law of the nonuniform
 L2-$1_{\sigma}$ scheme for time fractional Allen-Cahn model}
\author{
Hong-lin Liao\thanks{ORCID 0000-0003-0777-6832; School of Mathematics,
Nanjing University of Aeronautics and Astronautics,
Nanjing 211106, P.R. China; MIIT Key Laboratory of Mathematical Modelling and High 
Performance Computing of Air Vehicles, Nanjing 211106, P.R. China. Emails: liaohl@csrc.ac.cn and liaohl@nuaa.edu.cn.}
\quad Xiaohan Zhu\thanks{School of Mathematics, Nanjing University of Aeronautics and Astronautics,
211106, P.R. China. Email: cyzhuxiaohan@163.com.}
\quad Hong Sun\thanks{Department of Mathematics and Physics, 
	Nanjing Institute of Technology, Nanjing 211167; School of Mathematics, 
	Southeast University, Nanjing 210096, China. Email: sunhongzhal@126.com.}
}
\date{}
\maketitle
\normalsize

\begin{abstract}
We build an asymptotically compatible energy of 
the variable-step L2-$1_{\sigma}$ scheme 
for the time-fractional Allen-Cahn model 
with the Caputo's fractional derivative of order $\alpha\in(0,1)$, 
under a weak step-ratio constraint $\tau_k/\tau_{k-1}\geq r_{\star}(\alpha)$ for $k\ge2$, 	
where $\tau_k$ is the $k$-th time-step size and 
$r_{\star}(\alpha)\in(0.3865,0.4037)$ for $\alpha\in(0,1)$.
It provides a positive answer to the open problem in [J. Comput. Phys., 414:109473],
and, to the best of our knowledge, it is the first second-order nonuniform time-stepping scheme to preserve 
 both the maximum bound principle and the energy dissipation law of time-fractional Allen-Cahn model. 
The compatible discrete energy is constructed via a novel discrete gradient structure 
of the second-order L2-$1_{\sigma}$ formula by a local-nonlocal splitting technique.
It splits the discrete fractional derivative into two parts: one is a local term analogue to 
the trapezoid rule of the first derivative and the other is 
a nonlocal summation analogue to the L1 formula of  Caputo derivative.
Numerical examples with an adaptive time-stepping strategy are provided to 
show the effectiveness of our scheme and the asymptotic properties of the associated modified energy.

\noindent{\bf AMS subject classiffications.}\;\; 65M12, 65M06, 35Q99,   74A50
\end{abstract}

\section{Introduction}
\setcounter{equation}{0}

This paper continues to discuss the second-order nonuniform L2-$1_{\sigma}$ 
scheme in \cite{Liao2020ASecondorderJCP} for the time fractional Allen-Cahn (TFAC) model
\cite{DuYangZhou202TimeFractionalAllenCahn,HouXu2022sisc,LiWangYang:2017,
	LiuChengWangZhao2018TimefractionalACCH,
	Tang2018OnEnergyDissipation,ZhaoChenWang:2019},
\begin{align}\label{Cont: TFAC}
	\partial_t^\alpha \Phi= -\mu \quad\text{with}\quad
	\mu=\tfrac{\delta E}{\delta \Phi}=f(\Phi)-\epsilon^2\Delta\Phi,
\end{align}
where $\Phi$ is the phase variable  and $E$ is the Ginzburg-Landau type energy functional 
\begin{align}\label{Cont: energy functional}
E[\Phi]:=\int_{\Omega}\braB{\frac{\epsilon^2}{2}\abs{\nabla\Phi}^2
		+F(\Phi)}\zd\mathbf{x}
		\quad\text{with the potential $F(\Phi)=\frac{1}{4}\bra{\Phi^2-1}^2$}.
\end{align}
The variable $\Phi$ represents the concentration difference in a binary system on the domain 
$\Omega\subseteq\mathbb{R}^2$ and $\epsilon > 0$ is an interface width parameter. 
The notation $\partial_t^\alpha:={}_{0}^{C}\!D_{t}^{\alpha}$ represents
the Caputo's fractional  derivative of order $\alpha\in(0,1)$
with respect to $t$, defined by
\begin{align}\label{def:Caputo derivative}
	(\partial_{t}^{\alpha}v)(t)
	:=\int_{0}^{t}\omega_{1-\alpha}(t-s)v'(s)\zd{s}
	\quad \text{with $\omega_{\beta}(t):=t^{\beta-1}/\Gamma(\beta)$ for $\beta>0$.}
\end{align} 

As the fractional order $\alpha\rightarrow1^-$, the TFAC model 
\eqref{Cont: TFAC} is  asymptotically compatible  with the classical Allen-Cahn (AC) model \cite{AllenCahn1979AMicroscopicTheory, Feng2003NumericalAnalysisAllenCahn},
\begin{align}\label{Cont: classical Allen-Cahn}
	\partial_t \Phi=-\mu\quad\text{with}\quad
	\mu=\tfrac{\delta E}{\delta \Phi}.
\end{align}
As is well known, the AC model preserves the maximum bound principle 
\begin{align}\label{Cont: MaximumPrinciple}
	|\Phi(\mathbf{x},t)|\le{1}\;\text{for $t>0$}\quad
	\text{if}\quad|\Phi(\mathbf{x},0)|\le{1},
\end{align}
and the energy dissipation law
\begin{align}\label{Cont:AC energy dissipation law}
	\frac{\zd E}{\zd t}
	+\mynormb{\mu}^2=0\quad \text{for $t>0$},
\end{align}
where the $L^2$ norm $\mynorm{u}:=\sqrt{(u,u)}$ with 
the inner product $\bra{u,v}:=\int_{\Omega}uv\zd\mathbf{x}$
for $u,v\in L^2(\Omega).$

In recent years, the TFAC model \eqref{Cont: TFAC}
 has received extensive attentions and researches in 
accurately describing anomalous diffusion problems and developing 
efficient numerical methods to capture 
the long-time coarsening dynamics \cite{DuYangZhou202TimeFractionalAllenCahn,
	HouXu2022sisc,HouXu2022secondorderSAV,HouZhuXu2021L2SAVtfac,JiLiao2020SimpleMPTFAC,
	Karaa:2021,LiWangYang:2017,Liao2020ASecondorderJCP,LiuChengWangZhao2018TimefractionalACCH,
	QuanTangWang2022DecreasingUpperBound,Quan2020HowToDefine,
	Tang2018OnEnergyDissipation,ZhaoChenWang:2019}. 
A very interesting problem is whether the TFAC model \eqref{Cont: TFAC} also inherits   
the maximum bound principle \eqref{Cont: MaximumPrinciple} and 
the energy law \eqref{Cont:AC energy dissipation law} in an appropriate nonlocal form.
Tang, Yu and Zhou \cite{Tang2018OnEnergyDissipation} established the maximum bound principle
and derived a global energy dissipation law.
Quan, Tang and Yang proposed two nonlocal energy decaying laws for the time-fractional phase field models in \cite{Quan2020HowToDefine},
including a time-fractional energy dissipation law
$\bra{\partial_t^\alpha E}(t)\le 0$ for $t>0$,
and a weighted energy dissipation law,
$\frac{\zd E_{\omega}}{\zd t}\le 0$ for $t>0$
where $E_{\omega}(t):=\int_0^1\omega(s)E(s t)\zd s$ is a nonlocal weighted energy.
Du, Yang and Zhou \cite{DuYangZhou202TimeFractionalAllenCahn} 
investigated the well-posedness and the regularity of solution, and showed that
the solution satisfies $\mynorm{\partial_t^\alpha \Phi}_{L^p(0,T;L^2)}+\mynorm{\Delta\Phi}_{L^p(0,T;L^2)}\le C$ for any $p\in[2,2/\alpha)$ if the initial data $\Phi_0\in H^1_0(\Omega)$. They also discretized 
the fractional derivative by backward Euler convolution quadrature  and
developed some unconditionally solvable and stable time stepping schemes of order $O(\tau^\alpha)$, such as a convex
splitting scheme, a weighted convex splitting scheme and a linear weighted stabilized scheme.
A fractional energy dissipation law was also obtained in \cite{DuYangZhou202TimeFractionalAllenCahn}
and the discrete energy dissipation laws (in a weighted average sense) 
were derived for two weighted schemes on the uniform time mesh. 
Recently, Hou and Xu \cite{HouXu2022sisc,HouXu2022secondorderSAV,HouZhuXu2021L2SAVtfac} 
split the nonlocal time-fractional derivatives, including the L1-type, L2-$1_{\sigma}$ and L2 schemes, to local and nonlocal terms for time-fractional phase field equations, and treated
 the derived nonlocal term with the SAV technique,
so that the resulting modified discrete energy decays with respect to time. 


Naturally, we hope that the energy dissipation law of the TFAC model is asymptotically compatible with 
the energy dissipation law \eqref{Cont:AC energy dissipation law} of 
the AC equation \eqref{Cont: classical Allen-Cahn}, 
just as the TFAC model \eqref{Cont: TFAC} is asymptotically compatible 
with the AC equation \eqref{Cont: classical Allen-Cahn}. By reformulating 
the Caputo's form \eqref{Cont: TFAC} of the TFAC model  into the Riemann-Liouville 
form, one can follow the arguments \cite[Section 1]{LiaoTangZhou2021EnergyStableTFAC} to obtain the following
variational energy dissipation law
\begin{align}\label{Cont:TFAC variational energy dissipation law}
	\frac{\zd E_\alpha}{\zd t}
	+\frac{1}{2}\omega_{\alpha}(t)\mynormb{\mu}^2
	-\frac{1}{2} \int_{0}^{t}\omega_{\alpha-1}(t-s)\mynorm{\mu(t)-\mu(s)}^2\zd{s}=0
	\quad \text{for $t>0$,}
\end{align}
with the variational energy
\begin{align}\label{Cont:TFAC variational energy}
	E_\alpha[\Phi]:=E[\Phi]
	+\frac{1}{2}\int_{0}^{t}\omega_{\alpha}(t-s)\mynormb{\mu(s)}^2\zd{s}.
\end{align}
As remarked, the variational energy dissipation law
\eqref{Cont:TFAC variational energy dissipation law} is asymptotically compatible with 
the energy dissipation law \eqref{Cont:AC energy dissipation law} of the AC equation \eqref{Cont: classical Allen-Cahn}.
In recent works \cite{LiaoZhuWang:2021TFAC,JiZhuLiao2022TFCH,LiaoTangZhou2021EnergyStableTFAC}, the \emph{discrete gradient structures} of some variable-step L1-type formulas of
the Caputo and Riemann-Liouville fractional derivatives were 
constructed via the so-called discrete orthogonality convolution kernels. 
The  associated discrete energy dissipation laws of the corresponding schemes for time-fractional gradient flows were shown to be asymptotically compatible with the original energy law of their integer-order counterparts.
As seen, the variational energy \eqref{Cont:TFAC variational energy} and its discrete versions in \cite{LiaoZhuWang:2021TFAC,JiZhuLiao2022TFCH,LiaoTangZhou2021EnergyStableTFAC}
are always incompatible with the original (discrete) energy, in the fractional order limit $\alpha\rightarrow1^-$, due to the presence of the history effect in the Riemann-Liouville integral form.

In a recent interesting work \cite{QuanTangWang2022DecreasingUpperBound},
 Quan et al. proposed the following nonlocal-in-time modiﬁed energy (in our notations), formulated as the
 summation of the original energy $E[\Phi]$ and a new accumulation term due to the memory effect
 from time fractional derivative,
\begin{align}\label{Cont: QTW modified energy}
	\mathcal{E}_\alpha[\Phi]:=E[\Phi]
	+\frac{\omega_{1-\alpha}(t)}{2}\mynorm{\Phi(t)-\Phi(0)}^2
	-\frac{1}{2}\int_0^t\omega_{-\alpha}(t-s)\mynorm{\Phi(t)-\Phi(s)}^2\zd{s}.
\end{align}
It is remarkable that  this decreasing upper bound functional  
decays with respect to time and is asymptotically compatible with the original energy $E[\Phi]$
in the fractional order limit $\alpha\rightarrow1^-$.
The associated modified energy and energy dissipation laws at discrete time levels were derived in \cite{QuanTangWang2022DecreasingUpperBound} for 
several L1 and L2-type implicit-explicit stabilized schemes on the uniform time mesh in solving 
time-fractional gradient flow problems. The modified discrete energies of these time-stepping
schemes are incompatible (in the sense of $\alpha\rightarrow1^-$) 
with their integer-order counterparts, see 
\cite[Theorems 4.1-4.2 and 5.2]{QuanTangWang2022DecreasingUpperBound}, 
although they always coincide with the original energy at the initial time $t=0$ 
and as the time $t$ tends to $\infty$ (the steady state).

As observed from the extensive experiments in 
\cite{HouXu2022sisc,HouXu2022secondorderSAV,HouZhuXu2021L2SAVtfac,
	LiaoZhuWang:2021TFAC,JiZhuLiao2022TFCH,
	LiWangYang:2017,Liao2020ASecondorderJCP,LiaoTangZhou2021EnergyStableTFAC,
	LiuChengWangZhao2018TimefractionalACCH,
	QuanTangWang2022DecreasingUpperBound,Quan2020HowToDefine,
	Tang2018OnEnergyDissipation,ZhaoChenWang:2019},
the time-fractional gradient flow models admit 
multiple time scales in the long-time coarsening dynamics approaching the steady state. 
Thus an adaptive time-stepping strategy seems more suitable to
capture the multiscale behaviors and to save computational cost 
by choosing proper time-step sizes at different time periods. 
It is natural to require theoretically reliable time-stepping methods 
on general setting of time-step variations. 
We consider a time mesh $0=t_{0}<t_{1}<\cdots<t_{k}<\cdots<t_{N}=T$ 
with the time-step sizes $\tau_{k}:=t_{k}-t_{k-1}$ for $1\le{k}\le{N}$
and the time-step ratios $r_k:=\tau_k/\tau_{k-1}$ for $2\le k\le N$
(the step-ratio notation $\rho_k=1/r_{k+1}$ was used in previous works  \cite{Liao2021SecondorderReactionsudiffusion,Liao2020ASecondorderJCP}). 

The nonuniform L2-$1_{\sigma}$ formula \cite{Liao2021SecondorderReactionsudiffusion} 
of the Caputo derivative was applied in \cite{Liao2020ASecondorderJCP} to 
build a second-order maximum-bound-principle-preserving scheme for the TFAC model \eqref{Cont: TFAC};
however, no energy dissipation law of the second-order L2-$1_{\sigma}$ scheme has been established
due to the lack of technology to establish the positive definiteness (or \emph{discrete gradient structure}) of the 
discrete derivative on general nonuniform time meshes. 
Very recently, Quan and Wu \cite{QuanWu2022stabilityL21sigma} made a key progress 
on this issue by establishing the positive definiteness of 
the L2-$1_{\sigma}$ formula on general nonuniform meshes with 
a mild step-ratio restriction $r_k\geq 0.475329$ for $k\geq 2$. 
They also proved the $H^1$-stability of the L2-$1_{\sigma}$ 
time-stepping scheme for linear subdiffusion equations. 
Nonetheless, the positive definiteness of the discrete fractional operator 
should be inadequate to build an asymptotically compatible discrete 
energy for the TFAC model \eqref{Cont: TFAC}.

In this paper, we revisit the variable-step L2-$1_{\sigma}$ scheme for the TFAC model \eqref{Cont: TFAC}
and provide a positive answer to the unsolved problem raised in \cite{Liao2020ASecondorderJCP}.
Our contribution is two folds:
\begin{itemize}
	\item A new local-nonlocal splitting is adopted to split the L2-$1_{\sigma}$ formula into two parts:
	one is a local term analogue to the trapezoid rule of the first derivative and the other is 
	a nonlocal summation analogue to the L1 formula of  Caputo derivative. 
	Then a novel \emph{discrete gradient structure} (DGS)
	of the L2-$1_{\sigma}$ formula is constructed under an updated step-ratio constraint,
	that is $r_k\geq r_{\star}(\alpha)$ for $k\ge2$ 
	where $r_{\star}(\alpha)\in(0.3865,0.4037)$ for $\alpha\in(0,1)$, 
	which is a little weaker than those in previous studies 
	\cite{Liao2021SecondorderReactionsudiffusion,Liao2020ASecondorderJCP,QuanWu2022stabilityL21sigma}.
	\item An asymptotically compatible energy and the associated energy dissipation law of the variable-step L2-$1_{\sigma}$ scheme are established for the TFAC model \eqref{Cont: TFAC}. 
	To the best of our knowledge, it is the first second-order nonuniform time-stepping scheme to preserve 
	both the maximum bound principle and an asymptotically compatible energy dissipation law of the TFAC model. 
\end{itemize}

This paper is organized as follows. Section 2 describes the L2-$1_{\sigma}$ formula 
 and constructs a new discrete gradient structure by using a novel local-nonlocal splitting.
An asymptotically compatible discrete energy and the associated energy law of the L2-$1_{\sigma}$
implicit time-stepping scheme
are established in Section 3 for the TFAC model. 
Some numerical experiments are included in the last section to support our theoretical results.

\section{Discrete gradient structure of L2-1$_{\sigma}$ formula}\label{sec: Alikhanov}
\setcounter{equation}{0}

Set an off-set parameter $\theta :={\alpha}/{2}$ 
and $t_{k-\theta}:=\theta t_{k-1}+(1-\theta) t_{k}$.
For any grid function $\{w^{k}\}_{k=0}^N$, denote $\diff w^{k}:=w^{k}-w^{k-1}$, 
$\partial_{\tau}w^{k-\frac12}:=\diff w^{k}/\tau_k$ 
and the operator $w^{k-\theta}:=\theta w^{k-1}+(1-\theta) w^{k}$ for $k\geq 1$.
 For simplicity of presentation, we will always use the following notations
\[
\varpi_n(t):=-\omega_{2-\alpha}(t_{n-\theta}-t)
\quad\text{and} \quad
\varpi_n^{\prime}(t):=\omega_{1-\alpha}(t_{n-\theta}-t)
\quad\text{for $0\leq t\leq t_{n-\theta}$}.
\]

\subsection{Variable-step L2-1$_{\sigma}$ formula}

Let $\Pi_{1,k}v$ be the linear interpolant of a function~$v$ 
with respect to the two nodes $t_{k-1}$~and $t_k$,  and let $\Pi_{2,k}v$ denote the
quadratic interpolant with respect to $t_{k-1}$, $t_k$~and $t_{k+1}$. 
One has
\[
\bra{\Pi_{1,k}v}'(t)=\frac{\diff v^k}{\tau_k}
\;\;\text{and}\;\;
\bra{\Pi_{2,k}v}'(t)=\frac{\diff v^k}{\tau_k}
+\frac{2(t-t_{k-1/2})}{\tau_{k+1}(\tau_k+\tau_{k+1})}
\brab{\diff v^{k+1}-r_{k+1}\diff v^k}.
\]
The L2-1$_{\sigma}$ formula 
\cite{Alikhanov2015NewScheme,Liao2021SecondorderReactionsudiffusion,
	Liao2020ASecondorderJCP,QuanWu2022stabilityL21sigma} 
of  Caputo derivative~\eqref{def:Caputo derivative} can be obtained by
employing a quadratic interpolant in each subinterval $[t_{k-1}, t_k]$ for $1\leq k \leq n-1$
 and a linear interpolant in the final subinterval $[t_{n-1}, t_{n-\theta}]$.
That is
\begin{align}\label{eq: L2-1_sigma}
(\partial^{\alpha}_{\tau}v)^{n-\theta}
	&\defeq\int_{t_{n-1}}^{t_{n-\theta}}\varpi_n'(s)
		\bra{\Pi_{1,n}v}'(s)\zd{s}	+\sum_{k=1}^{n-1}\int_{t_{k-1}}^{t_k}
		\varpi_n'(s)\bra{\Pi_{2,k}v}'(s)\zd{s}\\
   &=a^{(n)}_0\diff v^n+\sum_{k=1}^{n-1}\braB{
	a^{(n)}_{n-k}\diff v^k
	+ \frac{\diff v^{k+1}-r_{k+1}\diff v^k}{r_{k+1}(1+r_{k+1})}\zeta^{(n)}_{n-k}},\nonumber
\end{align}
where $a_{n-k}^{(n)}$ and $\zeta_{n-k}^{(n)}$ are positive coefficients
 defined by
\begin{align}
&a^{(n)}_{n-k}\defeq\frac{1}{\tau_k}\int_{t_{k-1}}^{\min\{t_k,t_{n-\theta}\}}\varpi_n'(s)\zd{s}
\quad \text{for $1\le k\le n$},\label{eq: an}\\
&\zeta^{(n)}_{n-k}\defeq\frac{2}{\tau_k^2}\int_{t_{k-1}}^{t_k}
	(s-t_{k-\frac12})\varpi_n'(s)\zd{s}\quad \text{for $1\le k\le n$}. \label{eq: bn}
\end{align}

\subsection{Discrete gradient structure}

Notice that if $\alpha\rightarrow1^-$, then $\omega_{2-\alpha}(t)\to1$
whereas $\omega_{1-\alpha}(t)\to0$, uniformly for~$t$ in any compact
subinterval of the open half-line~$(0,\infty)$.  Thus,
$a^{(n)}_0=\omega_{2-\alpha}((1-\theta)\tau_n)/\tau_n\to1/\tau_n$
whereas $a^{(n)}_{n-k}\to0$~and $\zeta^{(n)}_{n-k}\to0$ for $1\le k\le n-1$.
It follows that 
$$(\dfd{\alpha}v)^{n-\theta}\;\;\to\;\; \partial_{\tau}v^{n-\frac12}\quad\text{as $\alpha\rightarrow 1^-$,} $$
so that the  L2-1$_{\sigma}$ formula~\eqref{eq: L2-1_sigma} tends to the
Crank-Nicolson approximation at the offset point $t_{n-1/2}$ of the first time derivative.
This asymptotic property inspires us to split the discrete fractional derivative 
\eqref{eq: L2-1_sigma} into two parts,
\begin{align}\label{eq: Alikhanov-local nonlocal splitting}
	(\partial^{\alpha}_{\tau}v)^{n-\theta}
	\triangleq&\,\underbrace{\frac{\alpha}{2-\alpha}a_0^{(n)}\diff v^{n}}
	+\underbrace{\sum_{k=1}^{n}\hat{a}_{n-k}^{(n)}\diff v^k}\quad \text{for $n\ge1$},\\
	&\,\hspace{1cm}J_{CN}^n\hspace{2cm} J_{L1}^n\nonumber
\end{align}
where the discrete convolution kernels $\hat{a}_{n-k}^{(n)}$ are defined by, $\hat{a}_0^{(1)}\defeq \frac{2(1-\alpha)}{2-\alpha}a_0^{(1)}$ and
\begin{equation}\label{eq: hat_a kernels}
\hat{a}^{(n)}_{n-k}\defeq\begin{cases}
	\frac{2(1-\alpha)}{2-\alpha}a^{(n)}_0+\frac{1}{r_n(1+r_n)}\zeta^{(n)}_1,
	&\text{for $k=n$, $n\geq 2$}\\[1ex]
	a^{(n)}_{n-k}+\frac{1}{r_k(1+r_k)}\zeta^{(n)}_{n-k+1}-\frac{1}{1+r_{k+1}}\zeta^{(n)}_{n-k},
	&\text{for $2\le k\le n-1$, $n\geq 3$}\\[1ex]
	a^{(n)}_{n-1}-\frac{1}{1+r_2}\zeta^{(n)}_{n-1},
	&\text{for $k=1$, $n\geq 2$.}
\end{cases}
\end{equation}

It is seen that
the local part $J_{CN}^n$ is similar to the Crank-Nicolson approximation of first time derivative 
and the second part $J_{L1}^n$ represents an L1-type formula of Caputo derivative 
\eqref{def:Caputo derivative}.
In general, the coefficient of $J_{CN}^n$ should be properly large 
so that $J_{CN}^n\rightarrow \partial_{\tau}v^{n-\frac12}$  and $J_{L1}^n\rightarrow0$  
as the fractional order $\alpha\rightarrow 1^-$, and it should also be properly small so that 
the remaining part $J_{L1}^n$ admits a discrete gradient structure under a weak step-ratio constraint. 

We will apply the splitting formulation \eqref{eq: Alikhanov-local nonlocal splitting} to 
derive a new discrete gradient structure of the L2-1$_\sigma$ formula \eqref{eq: L2-1_sigma},  
which plays an important role in establishing 
the discrete energy dissipation law of our numerical scheme. As we will see in 
Lemma \ref{lemma: decreasing and convexity} and Theorem \ref{thm: DGS},
the current choice $\frac{\alpha}{2-\alpha}a_0^{(n)}$ in the local part $J_{CN}^n$ would be subtle and 
key to obtain the minimum allowable lower bound $r_{\star}(\alpha)$ of the time-step ratios
for an asymptotically compatible discrete energy of 
the corresponding L2-1$_{\sigma}$ time-stepping scheme.

\begin{lemma}\label{lemma: minimum step ratios}
	Let $r_{\star}=r_{\star}(\alpha)$ 
	be the unique positive root of the equation
	$$h(r_{\star},\alpha):=2\sqrt{\frac{2(1-\alpha/2)r_{\star}}{1+\alpha+(1-\alpha/2)r_{\star}}
		+\frac{r_{\star}}{1+r_{\star}}}+3-\frac{1}{r_{\star}^2(1+r_{\star})}=0\quad\text{for $\alpha\in(0,1)$.}$$ 
It holds that $0.3865\approx r_{\star}(0)< r_{\star}(\alpha)< r_{\star}(1)\approx0.4037$\; for $\alpha\in(0,1)$.
	\end{lemma}
\begin{proof}
	It is not difficult to find that
	\begin{align*}
		\partial_{r_{\star}} h=&\,
		\braB{\frac{2r_{\star}}{\frac{2+2\alpha}{2-\alpha}+r_{\star}}+\frac{r_{\star}}{1+r_{\star}}}^{-\frac12}
		\kbra{\frac{4}{(\frac{2+2\alpha}{2-\alpha}+r_{\star})^2}\frac{1+\alpha}{2-\alpha}
			+\frac{1}{(1+r_{\star})^2}}+\frac{3r_{\star}+2}{r_{\star}^3(1+r_{\star})^2}>0,\\
		\partial_{\alpha} h=&\,
		-\frac{12r_{\star}}{(2-\alpha)^{2}(\frac{2+2\alpha}{2-\alpha}+r_{\star})^2}
		\braB{\frac{2r_{\star}}{\frac{2+2\alpha}{2-\alpha}+r_{\star}}+\frac{r_{\star}}{1+r_{\star}}}^{-\frac12}
		<0\quad\text{for $\alpha\in(0,1)$ and $r_{\star}>0$.}
	\end{align*}
	Due to the fact $h(1/4,\alpha)h(1/2,\alpha)<0$ for any $\alpha\in(0,1)$, 
	the equation $h(r_{\star},\alpha)=0$ has a unique positive root $r_{\star}=r_{\star}(\alpha)\in(1/4, 1/2)$. 
	Thanks to the implicit function theorem, 
	one has $r_{\star}'(\alpha)>0$	and the root $r_{\star}(\alpha)$ is increasing with respect to $\alpha\in(0,1)$.
	That is,  $r_{\star}(0)< r_{\star}(\alpha)< r_{\star}(1)$. Then we complete the proof by 
	solving the two equations $h(r_{\star}(0),0)=0$ and  $h(r_{\star}(1),1)=0$ numerically
	and find that $r_{\star}(0)\approx 0.3865$ and $r_{\star}(1)\approx 0.4037$.
\end{proof}

\begin{lemma}\label{lemma: decreasing and convexity}
	For the discrete kernels $\hat{a}_{n-k}^{(n)}$ in \eqref{eq: hat_a kernels}, 
	define the auxiliary convolution kernels 
	\begin{align}\label{eq: Alikhanov modified kernels A}
		A_{0}^{(n)}:=2\hat{a}_{0}^{(n)}\quad\text{and}\quad 
		A_{n-k}^{(n)}:=\hat{a}_{n-k}^{(n)}\quad\text{for $1\le k\le n-1$.}
	\end{align} 
For $r_{\star}=r_{\star}(\alpha)$ defined in Lemma \ref{lemma: minimum step ratios},
assume that the adjacent time-step ratios $r_k$ fulfill 
	\begin{align}\label{condition: step ratio restriction}
		r_k \geq r_{\star}(\alpha) \quad \text{for $k\geq 2$.}
	\end{align}	
Then the auxiliary convolution kernels $A_{n-k}^{(n)}$ satisfy
	\begin{itemize}
		\item[(a)]
		$ A^{(n)}_{n-k-1}> A^{(n)}_{n-k}>0$\; for $1\le k\le n-1$ $(n\ge2)$;	
		\item[(b)]
		$ A^{(n-1)}_{n-1-k}>A_{n-k}^{(n)}$\; for $1\le k\le n-1$ $(n\ge2)$;		
		\item[(c)]		$ A^{(n-1)}_{n-2-k}-A^{(n-1)}_{n-1-k}
		> A^{(n)}_{n-k-1}-A^{(n)}_{n-k}$\; for $1\le k\le n-2$ $(n\ge3)$.
	\end{itemize}
\end{lemma}
\begin{proof}See the lengthy but technical proof in  Section \ref{sec: kernel A property}.
	\end{proof}

\begin{theorem}\label{thm: DGS}
	Under the step-ratio constraint \eqref{condition: step ratio restriction},
 it holds that
	\begin{align*}
		2(\diff v^n)(\partial^\alpha_{\tau}v)^{n-\theta}
		= \mathcal{G}[\diff v^n]-\mathcal{G}[\diff v^{n-1}]
		+\mathcal{R}[\diff v^n]+\frac{2\alpha a_0^{(n)}}{2-\alpha}(\diff v^n)^2\quad \text{for $n\geq 1$,}
	\end{align*}
where the nonnegative functionals $\mathcal{G}$ and $\mathcal{R}$ are defined by
	\begin{align*}
	\mathcal{G}[w^n]&:=\sum_{j=1}^{n-1}\brab{A^{(n)}_{n-j-1}-A^{(n)}_{n-j}}
	\braB{\sum_{\ell=j+1}^nw^\ell}^2+A^{(n)}_{n-1}\braB{\sum_{\ell=1}^nw^\ell}^2\quad \text{for $n\geq 0$,}\\
	\mathcal{R}[w^n]&:= \sum_{j=1}^{n-1}
	\brab{A_{n-j-2}^{(n-1)}\!-A_{n-j-1}^{(n-1)}\!-A_{n-j-1}^{(n)}+A_{n-j}^{(n)}}\braB{\sum_{\ell=j+1}^{n-1}w^\ell}^2
	+\brab{A_{n-2}^{(n-1)}-A_{n-1}^{(n)}}\braB{\sum_{\ell=1}^{n-1}w^\ell}^2.
	\end{align*}
\end{theorem}
\begin{proof}Obviously, Lemma \ref{lemma: decreasing and convexity}  guarantees that the two functionals 
	$\mathcal{G}$ and $\mathcal{R}$ are nonnegative.
The kernel-splitting formula \eqref{eq: Alikhanov-local nonlocal splitting} gives $(\partial^\alpha_{\tau}v)^{n-\theta}=J_{CN}^n +J_{L1}^n.$
Obviously, 
the local part $J_{CN}^n $ can be handled easily, that is, 
$$2(\diff v^n)J_{CN}^n=\frac{2\alpha}{2-\alpha}a_0^{(n)}(\diff v^n)^2.$$ 
By following the proof of \cite[Lemma 2.3]{LiaoLiuZhao:2022FBDF2} with $\sigma_{\min}=0$,
it is not difficult to obtain the following equality for the nonlocal term $J_{L1}^n$, 
\begin{align*}
	2(\diff v^n)\sum_{k=1}^n&\hat{a}_{n-k}^{(n)}(\diff v^k)
	=\sum_{j=1}^{n-1}\brab{A^{(n)}_{n-j-1}-A^{(n)}_{n-j}}
		\braB{\sum_{\ell=j+1}^n\diff v^\ell}^2+A^{(n)}_{n-1}\braB{\sum_{\ell=1}^n\diff v^\ell}^2\\
	&-\sum_{j=1}^{n-2}\brab{A^{(n-1)}_{n-j-2}-A^{(n-1)}_{n-1-j}}
		\braB{\sum_{\ell=j+1}^{n-1}\diff v^\ell}^2
		-A^{(n-1)}_{n-2}\braB{\sum_{\ell=1}^{n-1}\diff v^\ell}^2+\mathcal{R}[\diff v^n].
\end{align*}
Then the desired DGS follows immediately.	
\end{proof}
\begin{remark}\label{remark: lower bound of step ratios}
	The new condition \eqref{condition: step ratio restriction}
	improves the lower bound of step ratios, $r_k\ge0.4753$, in the recent works
	on the stability of variable-step second-order formulas, see more details in
	 \cite[Theorem 3.2]{QuanWu2022stabilityL21sigma}, \cite[Theorem 2.1]{LiaoLiuZhao:2022FBDF2}
	 and \cite[Theorem 3.2 and Corollary 3.3]{QuanWu:2022H1stabilityL2}. 
	 As mentioned before, the properly small part $J_{CN}^n=\frac{\alpha}{2-\alpha}a_0^{(n)}\diff v^n$  
	 in the local-nonlocal splitting \eqref{eq: Alikhanov-local nonlocal splitting} is key  
	to achieve the minimum lower bound $r_{\star}(\alpha)$. 
	Can the variable-step fractional BDF2 formula \cite{LiaoLiuZhao:2022FBDF2}
	have an appropriate DGS under the new lower bound $r_{\star}(\alpha)$ of the step ratios? 
	It is an interesting issue worthy of further investigations although it is out of our current scope.
	Note that, the previous choice $\frac{1}{2-\alpha}a_0^{(n)}$ 
	 in \cite[(2.1)]{LiaoLiuZhao:2022FBDF2} for the local part of the fractional BDF2 formula 
	 is critical to obtain the maximum allowable upper bound $r^*(\alpha)$ of step ratios, 
	 see \cite[Lemmas 2.1-2.2]{LiaoLiuZhao:2022FBDF2}.
\end{remark}


\section{The L2-1$_\sigma$ scheme and energy dissipation law}
\setcounter{equation}{0}

Consider finite differences in spatial directions.
Let the uniform length $h:=L/M$ for some integer $M$ and the discrete grid
$\bar{\Omega}_h:=\big\{\mathbf{x}_{h}=(ih,jh)\,|\,0\le i,j\le M\big\}$.
Let $\Delta_h$ be the second-order difference approximation of Laplacian operator $\Delta$.
Let $\myinner{\cdot, \cdot} $ be the discrete inner product and $\mynorm{\cdot}$ 
denote the discrete $L^2$ norm.
Also, the maximum norm $\mynorm{v}_{\infty}:= \max_{\mathbf{x}_h\in \bar{\Omega}_h } \abs{v_h}$.

We revisit the following nonuniform L2-1$_\sigma$ scheme for the TFAC model \eqref{Cont: TFAC},
\begin{align}\label{scheme: implicit Alikhanov TFAC}
	\bra{\partial_\tau^\alpha \phi}^{n-\theta}=- \mu^{n-\theta}\quad\text{with}\quad
	\mu^{n-\theta}=f(\phi)^{n-\theta}-\epsilon^2\Delta_h \phi^{n-\theta} \quad\text{for $n\ge 1$,}
\end{align}
where the notation $f(\phi)^{n-\theta}\defeq \theta f(\phi^{n-1})+(1-\theta)f(\phi^{n})$.
The following result shows that the nonlinear scheme \eqref{scheme: implicit Alikhanov TFAC} is uniquely solvable 
and preserves the  maximum bound principle.
\begin{theorem}\label{thm: maximum bound principle}
Under the step-ratio condition \eqref{condition: step ratio restriction} with the step-size restriction
\begin{align}\label{condition: unique solvable and maximum principle}
\tau_n \leq \min\bigg\{\sqrt[\alpha]{\frac{\theta\omega_{2-\alpha}(1-\theta)}{2(1-\theta)}}, \sqrt[\alpha]{\frac{h^2\omega_{2-\alpha}(1-\theta)}{4\epsilon^2}}\bigg\},
\end{align}
the second-order nonuniform L2-1$_\sigma$ scheme \eqref{scheme: implicit Alikhanov TFAC} 
is uniquely solvable and preserves the maximum bound principle, that is, 
$\mynorm{\phi^n}_{\infty}\leq 1$ for $1\leq n\leq N$ 
if $\mynorm{\phi^0}_{\infty}\leq 1$. 
\end{theorem}
\begin{proof}\cite[Lemma 3.3]{Liao2020ASecondorderJCP} 
	gives the unique solvability of \eqref{scheme: implicit Alikhanov TFAC}. 
Also, the maximum bound principle is verified in \cite[Theorem 3.1]{Liao2020ASecondorderJCP}
under the step-ratio constraint $ r_k \geq 4/7$.  
The key point is the kernel monotonicity of the L2-1$_\sigma$  formula \eqref{eq: L2-1_sigma},
see \cite[(3.10) and (3.15)]{Liao2020ASecondorderJCP}. In our notations, it requires  $\hat{a}^{(n)}_{n-k-1}\ge\hat{a}^{(n)}_{n-k}$ for $1\le k\le n-2$,
which follows directly from the definition 
\eqref{eq: Alikhanov modified kernels A} and Lemma \ref{lemma: decreasing and convexity} (a). 
Thus the maximum bound principle also holds 
under the condition \eqref{condition: step ratio restriction}.
\end{proof}

As the fractional order $\alpha\to 1^{-}$, the  nonuniform L2-1$_\sigma$ scheme \eqref{scheme: implicit Alikhanov TFAC} degenerates into the classical Crank-Nicolson scheme for the AC model \eqref{Cont: classical Allen-Cahn}
\begin{align}\label{scheme: CN}
\partial_{\tau}\phi^{n-\frac12}=-\mu^{n-\frac12} \quad\text{with} \quad
\mu^{n-\frac12}=f(\phi)^{n-\frac12}-\epsilon^2\Delta_h\phi^{n-\frac12}.
\end{align}
It preserves the maximum bound principle if
$\tau_n \leq \min \{\frac12, \frac{h^2}{4\epsilon^2}\}$, 
see \cite[Theorem 1]{HouTangYang2017maximumprincipleCNTAC}.
Theorem \ref{thm: maximum bound principle} indicates that the time-step condition 
\eqref{condition: unique solvable and maximum principle} 
is asymptotically compatible as $\alpha\to 1^{-}$.

Now we define a new modified energy $\mathcal{E}_\alpha\kbra{\phi^{n}}$ as follows,
\begin{align}\label{def: modified energy}
	\mathcal{E}_\alpha[\phi^{n}]:=E[\phi^{n}]
	+\frac{1}{2}\myinner{\mathcal{G}[\diff\phi^{n}], 1}
	\quad\text{for $n\ge 1$,}
\end{align}
where the nonnegative functional $\mathcal{G}$ is defined in Theorem \ref{thm: DGS}
and $E\kbra{\phi^{n}}$ is 
the discrete counterpart of the Ginzburg-Landau energy functional \eqref{Cont: energy functional}, that is,
\begin{align*}
	E\kbra{\phi^{n}}:=\frac{\epsilon^2}{2}\mynormb{\nabla \phi^n}^2
	+\myinnerb{F(\phi^n),1}\quad\text{with}\quad
	F(\phi^n):=\frac14\brab{\brat{\phi^n}^2-1}^2\quad\text{for $n\ge 0$.}
\end{align*}
The following theorem says that the modified energy $\mathcal{E}_\alpha\kbra{\phi^{n}}$ decays at each time level.
\begin{theorem}\label{thm:Alikhanov energy law TFAC}
Under the step-ratio condition \eqref{condition: step ratio restriction} with the step-size restriction
\eqref{condition: unique solvable and maximum principle},
	the nonuniform L2-1$_\sigma$ scheme \eqref{scheme: implicit Alikhanov TFAC}
	preserves the following discrete energy dissipation law
	\begin{align*}
		\partial_\tau\mathcal{E}_\alpha\kbra{\phi^n}
		+\frac{\alpha}{2(2-\alpha)}a_0^{(n)}\tau_n\mynormb{\piff\phi^{n-\frac12}}^2\le 0 \quad \text{for $1\le n\le N$.}
	\end{align*}
\end{theorem}
\begin{proof}
Making the inner product of \eqref{scheme: implicit Alikhanov TFAC} by $\diff \phi^n$, one has 
\begin{align}\label{thmProof:Alikhanov energy law norm}
\myinnerb{(\partial^{\alpha}_{\tau}\phi)^{n-\theta}, \diff \phi^n}
+\epsilon^2\myinnerb{\nabla_h \phi^{n-\theta}, \nabla_h\diff \phi^n}
+\myinnerb{f(\phi)^{n-\theta}, \diff \phi^n}=0,
\end{align}
where the discrete Green's formula has been used. 
For the first term, Theorem \ref{thm: DGS} yields
	\begin{align*}
		\myinnerb{\bra{\partial_\tau^\alpha \phi}^{n-\theta},\diff \phi^n}
		\ge \frac12\myinner{\mathcal{G}[\diff\phi^{n}], 1}
		-\frac12\myinner{\mathcal{G}[\diff\phi^{n-1}], 1}
		+\frac{\alpha}{2-\alpha}a_0^{(n)}\mynormb{\diff \phi^n}^2.
	\end{align*}
With the fact $2a(a-b)=a^2-b^2+(a-b)^2$, the second term of 
\eqref{thmProof:Alikhanov energy law norm} can be formulated as
\begin{align*}
\epsilon^2\myinnerb{\nabla_h \phi^{n-\theta}, \nabla_h\diff \phi^n}=
\frac{\epsilon^2}{2}\mynorm{\nabla_h \phi^{n}}^2
-\frac{\epsilon^2}{2}\mynorm{\nabla_h \phi^{n-1}}^2
+\frac{(1-2\theta)\epsilon^2}{2}\mynorm{\nabla_h \diff\phi^{n}}^2.
\end{align*}
For any 
$ a,b \in [-1,1] $, we have the following two fundamental inequalities,
\begin{align*}
(a^3-a)(a-b)\geq\frac14[(a^2-1)^2-(b^2-1)^2]-(a-b)^2,\\
(b^3-b)(a-b)\geq\frac14[(a^2-1)^2-(b^2-1)^2]-(a-b)^2.
\end{align*}
With the help of Theorem \ref{thm: maximum bound principle}, 
the third term of \eqref{thmProof:Alikhanov energy law norm}  
can be bounded by
\begin{align*}
\myinnerb{f(\phi)^{n-\theta},\diff \phi^n}\geq 
\myinner{F(\phi^n),1}-\myinner{F(\phi^{n-1}),1}-\mynorm{\diff\phi^n}^2.
\end{align*}
	Inserting the above three estimates into \eqref{thmProof:Alikhanov energy law norm}, one gets
	\begin{align}\label{thmProof:Alikhanov energy law norm2}
	\braB{\frac{\alpha}{2(2-\alpha)}a_0^{(n)}-1}\mynormb{\diff\phi^n}^2
	+\frac{\alpha }{2(2-\alpha)}a_0^{(n)}\mynormb{\diff\phi^n}^2
		+\mathcal{E}_\alpha[\phi^{n}]
		\le \mathcal{E}_\alpha[\phi^{n-1}].
	\end{align}
The time-step condition \eqref{condition: unique solvable and maximum principle} implies that
 $\frac{\alpha}{2(2-\alpha)}a_0^{(n)}=\frac{\theta}{2(1-\theta)}\omega_{2-\alpha}(1-\theta)\tau_n^{-\alpha}\ge1$
 and the claimed result follows from \eqref{thmProof:Alikhanov energy law norm2} immediately.
\end{proof}

As the fractional index $\alpha\rightarrow 1^-$, it is easy to find that
$\hat{a}^{(n)}_{n-k}\to 0$ and $ {A}^{(n)}_{n-k}\to 0 $~for $1\le k\le n$.
The term $\mathcal{G}[\diff\phi^{n}] $ in the modified energy 
\eqref{def: modified energy} vanishes according to Theorem \ref{thm: DGS}.
Thus
\begin{align*}
\mathcal{E}_\alpha[\phi^{n}]=E[\phi^{n}]
	+\frac12\myinner{\mathcal{G}[\diff\phi^{n}], 1} \;\;\longrightarrow \;\; 
	E[\phi^{n}]\quad\text{as $\alpha\rightarrow 1^-$}
\end{align*}
such that the energy dissipation law in Theorem \ref{thm:Alikhanov energy law TFAC}
is asymptotically compatible with the energy dissipation 
law \cite[Theorem 2 or the inequality (3.21)]{HouTangYang2017maximumprincipleCNTAC}
of the Crank-Nicolson scheme \eqref{scheme: CN}, 
\begin{align*}
\partial_{\tau}\mathcal{E}_{\alpha}[\phi^n]+\frac{\alpha}{2(2-\alpha)}a_0^{(n)}\tau_n
\mynormb{\piff\phi^{n-\frac12}}^2\leq 0\;\;\longrightarrow \;\; 
\partial_{\tau}E[\phi^n]+\frac1{2}\mynormb{\mu^{n-\frac12}}^2\leq 0
\quad\text{as $\alpha\rightarrow 1^-$.}
\end{align*}
Obviously, the modified energy \eqref{def: modified energy}
of the  L2-1$_\sigma$ scheme also degrades into the original 
energy $E\kbra{\phi^n}$ in approaching the steady state ($\diff \phi^n\rightarrow0$), that is,
\[
\mathcal{E}_\alpha\kbra{\phi^n}\;\longrightarrow\;E\kbra{\phi^n}
\quad\text{as $t_n\rightarrow +\infty$.}
\]

\begin{remark}\label{remark: Alikanov error analysis}
	Although the nonuniform L2-1$_\sigma$ scheme \eqref{scheme: implicit Alikhanov TFAC} 
	proposed in \cite{Liao2020ASecondorderJCP} for the TFAC model is only considered in this paper; 
	however,  Theorem \ref{thm: DGS} would be useful 
	to derive the discrete energy dissipation laws of some other 
	L2-1$_\sigma$ approximations, such as convex splitting scheme
	and certain linearized schemes using the recent SAV techniques, 
	for time-fractional phase-field models.
   The interesting readers can refer to 
    \cite{Liao2018DiscreteGronwallInequality,KoptevaMeng:2020,Kopteva:2020,
    	Liao2021SecondorderReactionsudiffusion,QuanWu2022stabilityL21sigma}
     for the error analysis in the discrete $H^1$ or $L^2$ norm of 
    the variable-step L2-1$_\sigma$ time-stepping schemes. 
\end{remark}


\section{Proof of Lemma \ref{lemma: decreasing and convexity}}
\label{sec: kernel A property}
\setcounter{equation}{0}

This proof is rather technical and lengthy due to the nonuniform setting 
with the weak step-ratio condition \eqref{condition: step ratio restriction} 
and the inhomogeneity of the auxiliary convolution kernels $A_{n-k}^{(n)}$, which is defined by
\eqref{eq: Alikhanov modified kernels A} 
with the discrete kernels $\hat{a}_{n-k}^{(n)}$ in \eqref{eq: hat_a kernels}. 
To make the proof as clear as possible, we always avoid the usage of  $\hat{a}_{n-k}^{(n)}$ 
and directly use $a_{n-k}^{(n)}$ and $\zeta_{n-k}^{(n)}$ defined by \eqref{eq: an}-\eqref{eq: bn}.

Note that, the discrete convolution kernels of the nonuniform L2-1$_{\sigma}$ formula
\eqref{eq: L2-1_sigma} have been investigated in the previous work, 
see \cite[Theorem 2.1]{Liao2021SecondorderReactionsudiffusion}. However,
this proof is quite different from the proof of 
\cite[Theorem 2.1]{Liao2021SecondorderReactionsudiffusion} because 
the new step-ratio condition \eqref{condition: step ratio restriction} 
markedly improves the previous restriction $r_k\ge4/7\approx0.5714$
and the desired properties (b)-(c) are new. 
As the same point, we will continue to use the two integrals defined in
\cite[(4.1)]{Liao2021SecondorderReactionsudiffusion},
\begin{align*}
	I_{n-k}^{(n)}:=\int_{t_{k-1}}^{t_k}\frac{t_{k}-t}{\tau_k}\varpi_n''(t)\zd{t}
	\quad\text{and}\quad
	J_{n-k}^{(n)}:=\int_{t_{k-1}}^{t_k}\frac{t-t_{k-1}}{\tau_k}\varpi_n''(t)\zd{t}
	\quad\text{for $1\le k\le n-1$.}
\end{align*}
They will play a bridging role in building some useful links between the positive coefficients
$a_{n-k}^{(n)}$ and $\zeta_{n-k}^{(n)}$, see Lemmas \ref{lem: ank diff}-\ref{lem: I-J-theta qn2}.
It is worth mentioning that the new step-ratio condition \eqref{condition: step ratio restriction}
is a little better than those in the recent works \cite{LiaoLiuZhao:2022FBDF2,QuanWu2022stabilityL21sigma}.
Also, our proof seems more concise than the analysis in the proofs of 
\cite[Lemma 3.1 and Theorem 3.2]{QuanWu2022stabilityL21sigma,QuanWu:2022H1stabilityL2}.

\begin{lemma}\label{lem: ank diff}
	For $n\geq 2$, the positive coefficients $a^{(n)}_{n-k}$ in \eqref{eq: an} satisfy
	\begin{align*}
		&a_{n-k}^{(n)}-\varpi_n'(t_{k-1})=I_{n-k}^{(n)}\quad\text{for $1\leq k\leq n$,}\\
		&\varpi_n'(t_k)-a_{n-k}^{(n)}=J_{n-k}^{(n)}\quad\text{for $1\leq k\leq n-1$,}
	\end{align*}
	and then
	\begin{itemize}
		\item[(i)]
		$a_{n-k-1}^{(n)}-a_{n-k}^{(n)}
		=I_{n-k-1}^{(n)}+J_{n-k}^{(n)}$ \;for $1\le k\le n-1;$
		\item[(ii)]
		$\frac{4(1-\alpha)}{2-\alpha}a_{0}^{(n)}-a_{1}^{(n)}
		=\varpi_n'(t_{n-1})+J_{1}^{(n)}.$
	\end{itemize}
\end{lemma}
\begin{proof} This proof is similar to that of \cite[Lemma 4.6]{Liao2021SecondorderReactionsudiffusion}
	and we omit it here.	
\end{proof}
\begin{lemma}\label{lem: zeta-inequality}
	The positive coefficients $\zeta^{(n)}_{n-k}$ in \eqref{eq: bn} satisfy
	\begin{itemize}	
		\item[(i)] 
		$\zeta^{(n)}_{n-k}<\zeta^{(n-1)}_{n-k-1}$\; for $1\le k\le n-2;$
		\item[(ii)]
		$\zeta^{(n)}_{n-k-1}>r_{k+1}\zeta^{(n)}_{n-k}$\; for $1\le k\le n-2$, and
		$$\zeta^{(n)}_{n-k-1}- r_{k+1}\zeta^{(n)}_{n-k}<\zeta^{(n-1)}_{n-k-2}- r_{k+1}\zeta^{(n-1)}_{n-k-1}
		\quad\text{for $1\le k\le n-3$.}$$
	\end{itemize}	
\end{lemma}
\begin{proof} This proof is similar to that of \cite[Lemma 4.2]{LiaoLiuZhao:2022FBDF2}
	and we omit it here.	
\end{proof}

\begin{lemma}\label{lem: I-J-theta qn2}
	For $n\ge2$, define a discrete sequence $\{\beta_{k}\}$ as follows,
	\begin{align}\label{def: sequence beta n}
	\beta_{k}:=\frac{2(1-\alpha/2)r_{k}}{1+\alpha+(1-\alpha/2)r_{k}}\quad\text{for $k\ge2$}\,.
	\end{align}
	The positive coefficients $I^{(n)}_{n-k}, J_{n-k}^{(n)}$ and $\zeta_{n-k}^{(n)}$ satisfy
	\begin{itemize}
		\item[(i)] 
		$ I_{n-k}^{(n)}> (1+\beta_{k+1})\zeta^{(n)}_{n-k}$
		\; for $ 1 \leq k \leq n-1$, 
		and
		$$ I_{n-k}^{(n)}-(1+\beta_{k+1})\zeta^{(n)}_{n-k}
		<I_{n-k-1}^{(n-1)}-(1+\beta_{k+1})\zeta^{(n-1)}_{n-k-1}\quad 
		\text{for~ $1 \leq k \leq n-2$};$$		
	\item[(ii)]
		$ J_{n-k}^{(n)}> 3\zeta^{(n)}_{n-k}$\; for $1 \leq k \leq n-1$, and
		$$J_{n-k}^{(n)}-3\zeta^{(n)}_{n-k}<J_{n-k-1}^{(n-1)}-3\zeta^{(n-1)}_{n-k-1}
		\quad\text{for~ $1 \leq k \leq n-2$};$$
		\item[(iii)]
		$ r_n\zeta_{1}^{(n)}< \frac{\alpha}{3(2-\alpha)}\varpi_n'(t_{n-1})$\; for $ n\ge 2$.
	\end{itemize}
\end{lemma}
\begin{proof}
	By the integration by parts, one derives from \eqref{eq: bn} that
	\begin{align}\label{eq: eta nk formula 2}
		\zeta_{n-k}^{(n)}=\frac{1}{\tau_k^2}
		\int_{t_{k-1}}^{t_k}(t-t_{k-1})(t_k-t)\varpi_n''(t)\zd{t}\quad\text{for $ 1\le k\le n-1$,}
	\end{align}
	such that
	\begin{align*}
		I_{n-k}^{(n)}-\zeta_{n-k}^{(n)}
		=\frac{1}{\tau_k^2}\int_{t_{k-1}}^{t_k}(t_k-t)^2\varpi_n''(t)\zd{t}>0\quad\text{for $ 1\le k\le n-1$.}
	\end{align*}
It says that $I_{n-k}^{(n)}>\zeta_{n-k}^{(n)}$ for $ 1\le k\le n-1$. To obtain a sharper bound, 
we will compare  $\zeta_{n-k}^{(n)}$ with the difference term $I_{n-k}^{(n)}-\zeta_{n-k}^{(n)}$ by the Cauchy differential mean-value theorem. To do so, introduce two auxiliary functions with respect to $z\in[0,1]$,
\begin{align*}
	\zeta_{n,k}(z)\defeq&\frac{1}{\tau_k^2}\int_{t_{k-1}}^{t_{k-1}+z\tau_k}
	(t-t_{k-1})(t_{k-1}+z\tau_k-t)\varpi_n''(t)\zd{t}\quad\text{for $ 1\leq k\leq n-1$,}\\
	\Psi_{n,k}(z)\defeq&\frac{1}{\tau_k^2}\int_{t_{k-1}}^{t_{k-1}+z\tau_k}(t_{k-1}+z\tau_k-t)^2\varpi_n''(t)\zd{t} \quad\text{for $ 1\leq k\leq n-1$,}
\end{align*}
such that $\zeta_{n,k}(1)=\zeta_{n-k}^{(n)}$ and $\Psi_{n,k}(1)=I_{n-k}^{(n)}-\zeta_{n-k}^{(n)}$ for $ 1\leq k\leq n-1$.
Simple calculations give their derivatives
	\begin{align*}
		\zeta_{n,k}'(z)=&\frac{1}{\tau_k}\int_{t_{k-1}}^{t_{k-1}+z\tau_k}
		(t-t_{k-1})\varpi_n''(t)\zd{t},\quad
		\zeta_{n,k}''(z)=z\tau_k\varpi_n''(t_{k-1}+z\tau_k),\\
		\zeta_{n,k}'''(z)=&\tau_k\varpi_n''(t_{k-1}+z\tau_k)+z\tau_k^2\varpi_n'''(t_{k-1}+z\tau_k),
	\end{align*}
and
\begin{align*}
	\Psi_{n,k}'(z)=&\frac{2}{\tau_k}\int_{t_{k-1}}^{t_{k-1}+z\tau_k}(t_{k-1}+z\tau_k-t)\varpi_n''(t)\zd{t},
	\quad \Psi_{n,k}''(z)=2\int_{t_{k-1}}^{t_{k-1}+z\tau_k}\varpi_n''(t)\zd{t},\\
	\Psi_{n,k}'''(z)=&2\tau_k\varpi_n''(t_{k-1}+z\tau_k).
\end{align*}
	It is obvious that $\zeta_{n,k}(0)=0$, $\zeta_{n,k}'(0)=0$ and $\zeta_{n,k}''(0)=0 $ for $1\le k\le n-1$.
	Also, one has $\Psi_{n,k}(0)=0$, $\Psi_{n,k}'(0)=0$ and $\Psi_{n,k}''(0)=0$ for $1\le k\le n-1$.
	Thanks to the Cauchy differential mean-value theorem, there exist some $z_{1k},z_{2k}, z_{3k}\in(0,1)$ such that
	\begin{align*}
		\frac{I_{n-k}^{(n)}-\zeta_{n-k}^{(n)}}{\zeta_{n-k}^{(n)}}
		&=\frac{\Psi_{n,k}'(z_{1k})}{\zeta_{n,k}'(z_{1k})}=\frac{\Psi_{n,k}''(z_{2k})}{\zeta_{n,k}''(z_{2k})}
		=\frac{\Psi_{n,k}'''(z_{3k})}{\zeta_{n,k}'''(z_{3k})}
		=\frac{2}{1+(1+\alpha)\frac{z_{3k}\tau_k}{t_{n-\theta}-t_{k-1}-z_{3k}\tau_k}}\\
		&>\frac{2}{1+(1+\alpha)\frac{\tau_k}{t_{n-\theta}-t_{k}}}\ge\beta_{k+1}\quad\text{for $ 1\leq k\leq n-1$,}
	\end{align*}
where the fact $(1+\alpha)\varpi_n''(t)=(t_{n-\theta}-t)\varpi_n'''(t)$ for $0< t\le t_{n-1}$
was used in the last equality. 
Also, the decreasing property, $\varpi_n''(t)<\varpi_{n-1}''(t)$ for $0< t\le t_{n-2}$, implies that
\begin{align*}
		I_{n-k}^{(n)}-(1+\beta_{k+1})\zeta^{(n)}_{n-k}&=\int_{t_{k-1}}^{t_{k}}
		\frac{t_{k}-t}{\tau_k}\braB{\frac{t_{k}-t}{\tau_k}-\beta_{k+1}\frac{t-t_{k-1}}{\tau_k}}\varpi_n''(t)\zd{t}\\
		&<\int_{t_{k-1}}^{t_{k}}
		\frac{t_{k}-t}{\tau_k}\braB{\frac{t_{k}-t}{\tau_k}-\beta_{k+1}\frac{t-t_{k-1}}{\tau_k}}\varpi_{n-1}''(t)\zd{t}\\
		&=I_{n-k-1}^{(n-1)}-(1+\beta_{k+1})\zeta^{(n-1)}_{n-k-1}\quad\text{for $ 1\leq k\leq n-2$.}
	\end{align*}
The results in (i) are verified.  Moreover, the definitions of $J_{n-k}^{(n)}$ and \eqref{eq: eta nk formula 2} yield
	\begin{align*}
		J_{n-k}^{(n)}-3\zeta_{n-k}^{(n)}=
		\frac{3}{\tau_k^2}\int_{t_{k-1}}^{t_k}(t-t_{k-1})\brat{t-t_{k-1}-2\tau_k/3}\varpi_n''(t)\zd{t}>0
	\end{align*}
	for $1\leq k\leq n-1$ because $\varpi_n''(t)$ is increasing with repsect to $t$
	and  $$\int_{t_{k-1}}^{t_k}(t-t_{k-1})\brat{t-t_{k-1}-2\tau_k/3}\zd{t}=\tau_k^3\int_0^1s(s-2/3)\zd{s}=0.$$ 
	In the similar way, 
	\begin{align*}
		J_{n-k}^{(n)}-3\zeta_{n-k}^{(n)}<&\,
		\frac{3}{\tau_k^2}\int_{t_{k-1}}^{t_k}(t-t_{k-1})\brat{t-t_{k-1}-2\tau_k/3}\varpi_{n-1}''(t)\zd{t}\\
		=&\,J_{n-k-1}^{(n-1)}-3\zeta^{(n-1)}_{n-k-1}\quad \text{for $1\leq k\leq n-2$}.
	\end{align*}
	They confirm the result (ii).
	By using the fact $\alpha \varpi_n'(t)=(t_{n-\theta}-t)\varpi_n''(t)$ for $0\le t\le t_{n-1}$, 
	we take the case $k=n-1$ of \eqref{eq: eta nk formula 2} to find that 
	\begin{align*}
		\zeta_{1}^{(n)}\leq\varpi_n''(t_{n-1})\tau_{n-1}
		\int_{0}^{1}s(1-s)\zd{s}=\frac{\alpha}{6(1-\theta)r_n}\varpi_n'(t_{n-1}).
	\end{align*}
	It arrives at (iii) and completes the proof.
\end{proof}

We are in a position to present a detail proof of Lemma \ref{lemma: decreasing and convexity}.


\begin{proof}This proof contains three parts: 
	\textbf{Part I} gives the monotonic property (a), \textbf{Part II} checks the algebraic convexity (c) and \textbf{Part III} verifies the monotonic property (b).

	\textbf{(Part I)}\quad	The definition \eqref{eq: hat_a kernels} of the discrete kernels $\hat{a}_{n-k}^{(n)}$,
	will be used frequently in this proof. By using the definitions \eqref{eq: an} and \eqref{eq: bn}, one has $A_0^{(n)}=2\hat{a}_{0}^{(n)}>0$ and 
\begin{align*}
	A_{n-k}^{(n)}>a_{n-k}^{(n)}-\frac{1}{1+r_{k+1}}\zeta_{n-k}^{(n)}
	=\int_{t_{k-1}}^{t_k}\frac{t_{k+1}+t_k-2t}{(1+r_{k+1})\tau_k^2}\varpi_n'(t)\zd{t}>0
	\quad\text{for $1\le k\le n-1$,}
\end{align*}
due to the fact $\int_{t_{k-1}}^{t_k}(t_{k+1}+t_k-2t) \zd{t}>0$.
Thus the discrete kernels $A_{n-k}^{(n)}>0$ for $1\le k\le n.$

To check the monotonic property of $A_{n-k}^{(n)}$ with respect to $k$, 
we need the following two inequalities. 
If $r_l\ge r_{\star}(\alpha)$ for $l\ge2$, 
one applies  the definition \eqref{def: sequence beta n} and Lemma \ref{lemma: minimum step ratios} to get
\begin{align}\label{ieq: hybid ratio condition1}
	&\braB{\beta_{k+2}+\frac{r_{k+2}}{1+r_{k+2}}}r_{k+1}
	+\frac{1}{r_{k+1}}+3-\frac{1}{r_k^2(1+r_k)}\nonumber\\
	&\hspace{1cm}\geq 2\sqrt{\frac{2(1-\alpha/2)r_{k+2}}{1+\alpha+(1-\alpha/2)r_{k+2}}+\frac{r_{k+2}}{1+r_{k+2}}}
	+3-\frac{1}{r_{\star}^2(1+r_{\star})}\nonumber\\
	&\hspace{1cm}\geq 2\sqrt{\frac{2(1-\alpha/2)r_{\star}}{1+\alpha+(1-\alpha/2)r_{\star}}+\frac{r_{\star}}{1+r_{\star}}}
	+3-\frac{1}{r_{\star}^2(1+r_{\star})}=h(r_{\star},\alpha)=0\quad\text{for $k\ge2$,}
\end{align}
and, similarly, 
\begin{align}\label{ieq: hybid ratio condition2}
	\frac{3(2-\alpha)}{\alpha}r_{k+1}&+\frac{2+r_{k+1}}{r_{k+1}(1+r_{k+1})}
	+3-\frac{1}{r_k^2(1+r_{k})}
	\nonumber\\	&\ge2\sqrt{\frac{3(2-\alpha)}{\alpha}}+3-\frac{1}{r_{\star}^2(1+r_{\star})}
	>h(r_{\star},\alpha)=0\quad\text{for $k\ge2$.}
\end{align}
We will verify $A_{n-k-1}^{(n)}\ge A_{n-k}^{(n)}$ $(n\ge k+1\ge2)$ 
by the four separate cases:
(I.1) $k=1$ for $n=2$, (I.2) $k=1$ for $n\geq3$, (I.3) $2\leq k\leq n-2$ for $n\geq4$,
and  (I.4) $k=n-1$ for	$n\geq3$. 
\begin{itemize}
	\item[\textbf{(I.1)}] \textbf{The case $k=1$  ($n=2$)}. Lemma \ref{lem: ank diff} (ii) gives
	\begin{align}\label{eq: property a k1n2}
		A^{(2)}_0-A^{(2)}_1&=\frac{4(1-\alpha)}{2-\alpha}a^{(2)}_0-a^{(2)}_1+\frac{2+r_2}{r_2(1+r_2)}\zeta_1^{(2)}
		=\varpi_2'(t_{1})+J_{1}^{(2)}+\frac{2+r_2}{r_2(1+r_2)}\zeta_1^{(2)}.
	\end{align}
The desired inequality $A^{(2)}_0>A^{(2)}_1$ is obvious.
\item[\textbf{(I.2)}] \textbf{The case $k=1$ ($n\ge3$)}. 
Applying Lemma \ref{lem: ank diff} (i), we reformulate the difference term
$A_{n-2}^{(n)}-A_{n-1}^{(n)}$ into a linear combination with nonnegative coefficients, 
\begin{align}\label{eq: property a k1n3}
	A^{(n)}_{n-2}-A^{(n)}_{n-1}
	=&\,a_{n-2}^{(n)}-a_{n-1}^{(n)}+\frac{\zeta_{n-1}^{(n)}}{r_2}-\frac{\zeta_{n-2}^{(n)}}{1+r_3}
	=I_{n-2}^{(n)}+J_{n-1}^{(n)}+\frac{\zeta_{n-1}^{(n)}}{r_2}-\frac{\zeta_{n-2}^{(n)}}{1+r_3}\nonumber\\
	=&\,\bra{I_{n-2}^{(n)}-(1+\beta_{3})\zeta_{n-2}^{(n)}}
	+\brab{\beta_{3}+\frac{r_{3}}{1+r_{3}}}\zeta_{n-2}^{(n)}\nonumber\\
	&\,	+\brab{J_{n-1}^{(n)}-3\zeta_{n-1}^{(n)}}+\frac{(1+3r_2)}{r_2}\zeta_{n-1}^{(n)}.
\end{align}
Then Lemma \ref{lem: I-J-theta qn2} (i)-(ii) yields $A_{n-2}^{(n)}>A_{n-1}^{(n)}$ immediately.	
\item[\textbf{(I.3)}] \textbf{The general cases $k=2,3,\cdots, n-2$  ($n\ge4$)}. 
Lemma \ref{lem: ank diff} (i) gives
\begin{align}\label{eq: property a kn-2n4}
	A_{n-k-1}^{(n)}-A_{n-k}^{(n)}=&\,
	a_{n-k-1}^{(n)}-a_{n-k}^{(n)}+\frac{\zeta^{(n)}_{n-k}}{r_{k+1}}
	-\frac{\zeta_{n-k-1}^{(n)}}{1+r_{k+2}}
	-\frac{\zeta^{(n)}_{n-k+1}}{r_{k}(1+r_{k})}\nonumber\\
	=&\,\braB{I_{n-k-1}^{(n)}-(1+\beta_{k+2})\zeta_{n-k-1}^{(n)}}
	+\brab{\beta_{k+2}+\frac{r_{k+2}}{1+r_{k+2}}}
	\brab{\zeta_{n-k-1}^{(n)}-r_{k+1}\zeta^{(n)}_{n-k}}	\nonumber\\
	&\,+\brab{J_{n-k}^{(n)}-3\zeta^{(n)}_{n-k}}
	+\frac{1}{r_{k}^2(1+r_{k})}\brab{\zeta^{(n)}_{n-k}-r_k\zeta^{(n)}_{n-k+1}}\nonumber\\
	&\,+\kbra{\brab{\beta_{k+2}+\frac{r_{k+2}}{1+r_{k+2}}}r_{k+1}
		+\frac{1}{r_{k+1}}+3-\frac{1}{r_{k}^2(1+r_{k})}}
	\zeta^{(n)}_{n-k}
\end{align}	
for $2\le k\le n-2$. With the help of Lemma \ref{lem: zeta-inequality} (ii) and 
Lemma \ref{lem: I-J-theta qn2} (i)-(ii) together with the inequality 
\eqref{ieq: hybid ratio condition1}, the nonnegative linear combination
\eqref{eq: property a kn-2n4} leads to the desired inequality.	
\item[\textbf{(I.4)}] \textbf{The case $k=n-1$  ($n\ge3$)}. 
By using Lemma \ref{lem: ank diff} (ii), we have
\begin{align}\label{eq: property a kn-1n3}
	A_{0}^{(n)}-A_{1}^{(n)}=&\,\frac{4(1-\alpha)}{2-\alpha}a^{(n)}_0-a_{1}^{(n)}
	+\frac{(2+r_n)\zeta^{(n)}_{1}}{r_n(1+r_n)}-\frac{\zeta^{(n)}_{2}}{r_{n-1}(1+r_{n-1})}\nonumber\\
	=&\, \braB{\varpi_n'(t_{n-1})-\frac{3(2-\alpha)}{\alpha}r_{n}\zeta^{(n)}_{1}}
	+\brab{J_{1}^{(n)}-3\zeta^{(n)}_{1}}
	+\frac{\zeta^{(n)}_{1}-r_{n-1}\zeta^{(n)}_{2}}{r_{n-1}^2(1+r_{n-1})}\nonumber\\
	&\,+\kbra{\frac{3(2-\alpha)}{\alpha}r_{n}+\frac{2+r_n}{r_{n}(1+r_{n})} +3-\frac{1}{r_{n-1}^2(1+r_{n-1})}}\zeta^{(n)}_{1}.
\end{align}
Thus Lemma \ref{lem: zeta-inequality} (ii) and 
Lemma \ref{lem: I-J-theta qn2} (ii)-(iii) together with the constraint 
\eqref{ieq: hybid ratio condition2} arrive at the desired inequality, 
$A_{0}^{(n)}>A_{1}^{(n)}$.	It completes the proof of Lemma \ref{lemma: decreasing and convexity} (a).	
\end{itemize}

\textbf{(Part II)}\quad	This part verifies Lemma \ref{lemma: decreasing and convexity} (c)
with the help of a positive auxiliary function
\begin{align}\label{def: auxiliary function g}
	g(x,\alpha):=1+\frac{2(1-\alpha)x}{2-\alpha}
	-\frac{\alpha}{2-\alpha}\braB{x+\frac{2}{2-\alpha}}^{-\alpha}>0
	\quad\text{for $x>0$ and $0<\alpha<1$,}
\end{align}
due to the fact $\partial_xg(x,\alpha)>0$ and then $g(x,\alpha)>g(0,\alpha)>0$.
To make full use of the above linear combination 
forms \eqref{eq: property a k1n2}-\eqref{eq: property a kn-1n3}
of the difference term $A^{(n)}_{n-k-1}-A^{(n)}_{n-k}$ for different indexes $k$,
the algebraic convexity (c) will be checked via the equivalent form,
\begin{align}\label{ieq: equivalent Theorem (c)}
	A^{(n)}_{n-k-1}-A^{(n)}_{n-k} <A^{(n-1)}_{n-2-k}-A^{(n-1)}_{n-1-k}
	\quad\text{for $1\le k\le n-2$ $(n\ge3)$}.
\end{align}
To cover all of possibilities, we will consider the following four separate cases:
(II.1) $k=1$ for $n=3$, (II.2) $k=1$ for $n\geq4$, 
(II.3) $2\leq k\leq n-3$ for $n\geq5$, and  (II.4) $k=n-2$ for $n\geq4$.

\begin{itemize}
	\item[\textbf{(II.1)}] \textbf{The case $k=1$ ($n=3$)}. 
	Applying Lemma \ref{lem: zeta-inequality} (i) and Lemma \ref{lem: I-J-theta qn2} (i)-(ii)
	to the last two terms of 
	the linear combination form \eqref{eq: property a k1n3} for the case $n=3$,
	but retaining and merging the first two terms, one has	
	\begin{align*}
		A_{1}^{(3)}-A_{2}^{(3)}
		<&\,I_{1}^{(3)}-\frac{\zeta_{1}^{(3)}}{1+r_{3}}		
			+\brab{J_{1}^{(2)}-3\zeta_{1}^{(2)}}+\frac{(1+3r_2)}{r_2}\zeta_{1}^{(2)}\,.
	\end{align*}
 Then subtracting the equality \eqref{eq: property a k1n2} 
from the above inequality, we apply the first equality (taking $k=2$ and $n=3$) of Lemma \ref{lem: ank diff} to get
\begin{align*}
	A_{1}^{(3)}-A_{2}^{(3)}-\brab{A^{(2)}_0-A^{(2)}_1}
	<&\,I_{1}^{(3)}-\frac{\zeta_{1}^{(3)}}{1+r_{3}}-\varpi_2'(t_{1})
	=a_{1}^{(3)}-\varpi_3'(t_{1})-\frac{\zeta_{1}^{(3)}}{1+r_{3}}-\varpi_2'(t_{1})\\
	=&\,-\frac{(1-\theta)^{-\alpha}\tau_{2}^{-\alpha}}{(1+r_3)\Gamma(1-\alpha)}
	\braB{g(r_3,\alpha)+\frac{\alpha r_{3}}{2-\alpha}}<0,	
\end{align*}
where $g$ is the positive function defined by \eqref{def: auxiliary function g}.
	\item[\textbf{(II.2)}] \textbf{The case $k=1$ ($n\ge4$)}. 
	Applying the linear combination form \eqref{eq: property a k1n3},
	Lemma \ref{lem: zeta-inequality} (i) and Lemma \ref{lem: I-J-theta qn2} (i)-(ii), we have
	\begin{align*}
		A_{n-2}^{(n)}-A_{n-1}^{(n)}
		<&\,\braB{I_{n-3}^{(n-1)}-(1+\beta_{3})\zeta_{n-3}^{(n-1)}}
		+\brab{\beta_{3}+\frac{r_{3}}{1+r_{3}}}\zeta_{n-3}^{(n-1)}\nonumber\\
		&\,+\brab{J_{n-2}^{(n-1)}-3\zeta_{n-2}^{(n-1)}}	
		+\frac{(1+3r_2)}{r_2}\zeta_{n-2}^{(n-1)}
		=A_{n-3}^{(n-1)}-A_{n-2}^{(n-1)}.
	\end{align*}
	\item[\textbf{(II.3)}] \textbf{The general cases $k=2,3,\cdots, n-3$  ($n\ge5$)}. By using the equality 
	\eqref{eq: property a kn-2n4}, Lemma \ref{lem: zeta-inequality} and 
	Lemma \ref{lem: I-J-theta qn2} (i)-(ii) together with the step-ratio 
	constraint \eqref{ieq: hybid ratio condition1}, we obtain that
	\begin{align}\label{eq: column property a kn-2n4}
		A_{n-k-1}^{(n)}-A_{n-k}^{(n)}<&\,
		\braB{I_{n-k-2}^{(n-1)}-(1+\beta_{k+2})\zeta_{n-k-2}^{(n-1)}}
		+\brab{\beta_{k+2}+\frac{r_{k+2}}{1+r_{k+2}}}
		\brab{\zeta_{n-k-2}^{(n-1)}-r_{k+1}\zeta^{(n-1)}_{n-k-1}}	\nonumber\\
		&\,+\brab{J_{n-k-1}^{(n-1)}-3\zeta^{(n-1)}_{n-k-1}}
		+\frac{1}{r_{k}^2(1+r_{k})}\brab{\zeta^{(n-1)}_{n-k-1}-r_k\zeta^{(n-1)}_{n-k}}\nonumber\\
		&\,+\kbra{\brab{\beta_{k+2}+\frac{r_{k+2}}{1+r_{k+2}}}r_{k+1}
			+\frac{1}{r_{k+1}}+3-\frac{1}{r_{k}^2(1+r_{k})}}\zeta^{(n-1)}_{n-k-1}\nonumber\\
		=&\,A_{n-k-2}^{(n-1)}-A_{n-k-1}^{(n-1)}\qquad\text{for $2\le k\le n-3$.}
	\end{align}	
	\item[\textbf{(II.4)}] \textbf{The case $k=n-2$  ($n\ge4$)}. Applying Lemma \ref{lem: zeta-inequality}, 
	Lemma \ref{lem: I-J-theta qn2} (ii) and the 
	restriction \eqref{ieq: hybid ratio condition1} to the last three terms of 
	the linear combination form \eqref{eq: property a kn-2n4} for the case $k=n-2$
	(but retaining and merging the first two terms), we obtain that
	\begin{align}\label{eq: part column property a kn-2n4}
		A_{1}^{(n)}-A_{2}^{(n)}<&\,I_{1}^{(n)}-\frac{\zeta_{1}^{(n)}}{1+r_{n}}
		-\brab{\beta_{n}+\frac{r_{n}}{1+r_{n}}}r_{n-1}\zeta^{(n)}_{2}	\nonumber\\
		&\,+\brab{J_{1}^{(n-1)}-3\zeta^{(n-1)}_{1}}
		+\frac{1}{r_{n-2}^2(1+r_{n-2})}\brab{\zeta^{(n-1)}_{1}-r_{n-2}\zeta^{(n-1)}_{2}}\nonumber\\
		&\,+\kbra{\brab{\beta_{n}+\frac{r_{n}}{1+r_{n}}}r_{n-1}
			+\frac{1}{r_{n-1}}+3-\frac{1}{r_{n-2}^2(1+r_{n-2})}}\zeta^{(n-1)}_{1}\,.
	\end{align}
	Replacing the index $n$ with $n-1$ in the formulation \eqref{eq: property a kn-1n3}, one has
	\begin{align}\label{eq: part column property a kn-1n3}
		A_{0}^{(n-1)}-A_{1}^{(n-1)}
		=&\, \varpi_{n-1}'(t_{n-2})+\brab{J_{1}^{(n-1)}-3\zeta^{(n-1)}_{1}}
		+\frac{\zeta^{(n-1)}_{1}-r_{n-2}\zeta^{(n-1)}_{2}}{r_{n-2}^2(1+r_{n-2})}\nonumber\\
		&\,+\kbra{\frac{2+r_{n-1}}{r_{n-1}(1+r_{n-1})} +3-\frac{1}{r_{n-2}^2(1+r_{n-2})}}\zeta^{(n-1)}_{1}.
	\end{align}
	 By subtracting \eqref{eq: part column property a kn-1n3} 
	from \eqref{eq: part column property a kn-2n4} and dropping some nonpositive terms, 
	we apply the first equality (taking $k=n-1$) of Lemma \ref{lem: ank diff}
	and Lemma \ref{lem: I-J-theta qn2} (iii) to derive that	
	\begin{align}\label{eq: algebraic convexity at ends}
		&\,A_{1}^{(n)}-A_{2}^{(n)}-\brab{A_{0}^{(n-1)}-A_{1}^{(n-1)}}
		<I_{1}^{(n)}-\frac{\zeta_{1}^{(n)}}{1+r_{n}}
		-\varpi_{n-1}'(t_{n-2})
		+\brab{\beta_{n}+\frac{r_{n}}{1+r_{n}}}r_{n-1}\zeta^{(n-1)}_{1}\nonumber\\
		&\,\hspace{2cm}<a_{1}^{(n)}-\varpi_n'(t_{n-2})-\frac{\zeta_{1}^{(n)}}{1+r_{n}}
		-\varpi_{n-1}'(t_{n-2})
		+\brab{\beta_{n}+\frac{r_{n}}{1+r_{n}}}\frac{\alpha\varpi_{n-1}'(t_{n-2})}{3(2-\alpha)}\nonumber\\
		&\,\hspace{2cm}<a_{1}^{(n)}-\varpi_n'(t_{n-2})-\frac{\zeta_{1}^{(n)}}{1+r_{n}}
		-\varpi_{n-1}'(t_{n-2})+\frac{r_{n}}{1+r_{n}}\frac{\alpha}{2-\alpha}\varpi_{n-1}'(t_{n-2})\nonumber\\
		&\,\hspace{2cm}=-\frac{(1-\theta)^{-\alpha}\tau_{n-1}^{-\alpha}}{(1+r_n)\Gamma(1-\alpha)}g(r_n,\alpha)<0,
	\end{align} 
	where the fact $\beta_{n}=\frac{2(1-\alpha/2)r_{n}}{1+\alpha+(1-\alpha/2)r_{n}}< \frac{2r_{n}}{1+r_{n}}$ 
	has been used in the third inequality and $g$ is the positive function  defined by \eqref{def: auxiliary function g}.
	It completes the proof of Lemma \ref{lemma: decreasing and convexity} (c).
\end{itemize}

\textbf{(Part III)}\quad It remains to prove Lemma \ref{lemma: decreasing and convexity} (b). 
The definition of $A_{n-1}^{(n)}$ gives
\begin{align*}
	A_{n-1}^{(n)}=a_{n-1}^{(n)}-\frac{1}{1+r_2}\zeta_{n-1}^{(n)}
	\le\frac{1}{\tau_1}\int_{t_0}^{t_1}\frac{t_2+t_1-2t}{(1+r_2)\tau_1}\varpi_{n-1}'(t) \zd{t}=A_{n-2}^{(n-1)}
	\quad\text{for $n\ge3$.}
\end{align*}
For the simple case $n=2,$ a direct calculation arrives at 
\begin{align*}
	A_{1}^{(2)}-A_{0}^{(1)}=&\,a_{1}^{(2)}-\frac{1}{1+r_2}\zeta_{1}^{(2)}-\frac{4(1-\alpha)}{2-\alpha}a_{0}^{(1)}
	=a_{1}^{(2)}-\frac{1}{1+r_2}\zeta_{1}^{(2)}-2\varpi_{1}'(t_{0})\\
	=&\,-\frac{(1-\theta)^{1-\alpha}\tau_{1}^{-\alpha}}{(1+r_2)\Gamma(1-\alpha)}
	\kbra{\frac{2\alpha}{2-\alpha}+2\braB{r_2+\frac{2}{2-\alpha}}-\braB{r_2+\frac{2}{2-\alpha}}^{1-\alpha}}<0.
\end{align*}
Thus one has $A_{n-2}^{(n-1)}>A_{{n-1}}^{(n)}$ for $n\ge2.$
Thanks to \eqref{ieq: equivalent Theorem (c)}, we have 
\begin{align*}
	A^{(n-1)}_{n-k-2}-A^{(n)}_{n-k-1}>A^{(n-1)}_{n-k-1}-A^{(n)}_{n-k}>0\quad\text{for $1\le k\le n-2$ $(n\ge3)$,}
\end{align*}
which leads to the monotonic property (b) directly.
The proof of Lemma \ref{lemma: decreasing and convexity} is completed.
\end{proof}

\section{Numerical examples}
\setcounter{equation}{0}

In this section, we present several numerical examples
to illustrate the efficiency and accuracy of the nonuniform
L2-1$_\sigma$ method \eqref{scheme: implicit Alikhanov TFAC} for the TFAC model \eqref{Cont: TFAC}.
At each time level, the nonlinear equation is solved by employing  a fixed-point algorithm with the termination error $10^{-12}$.
Also, the sum-of-exponentials technique\cite{Jiang2017FastEvaluationCaputo} with the absolute tolerance error $\varepsilon=10^{-12}$ and cut-off time $\Delta{t}=10^{-12}$ is used to speed up the convolution computation of the L2-1$_\sigma$ formula \eqref{eq: L2-1_sigma}.

\subsection{Accuracy verification}

\begin{example}\label{exam:accuracy test}
Consider an exact solution
$\Phi(\mathbf{x},t)=\omega_{1+\sigma}(t)\sin({ x})\sin({ y})$
with $\sigma\in(0,1)$ by adding an exterior force to the TFAC model \eqref{Cont: TFAC}.
\end{example}

We use a $512^2$ spatial mesh to discretize the spatial domain $(0,2\pi)^2$ and take $T=1$ and $\epsilon^2=0.1$.
The time interval $[0,T]$ is always divided into two parts $[0, T_0]$ and $[T_0, T]$ with total $N$ subintervals, where $T_0=\min\{1/\gamma,T\}$ and $N_0=\lceil \frac{N}{T+1-\gamma^{-1}}\rceil$.
The graded mesh $t_{k}=T_{0}(k/N_0)^{\gamma}$ is applied inside 
the initial part $[0,T_{0}]$ for resolving the initial singularity.
The random time-steps $\tau_{N_{0}+k}:=(T-T_{0})\epsilon_{k}/S_1$ for $1\le k\le N_1$ are used in the remainder interval $[T_{0},T]$ where $N_1:=N-N_0$, $S_1:=\sum_{k=1}^{N_1}\epsilon_{k}$ and $\epsilon_{k}\in(0,1)$ are the random numbers.

\begin{figure}[htb!]
	\centering
	\includegraphics[scale=0.4]{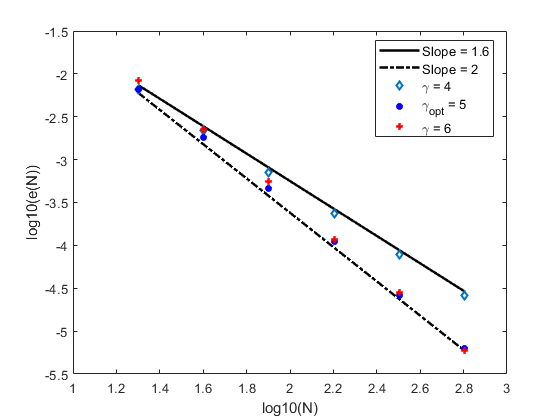} 
	\includegraphics[scale=0.4]{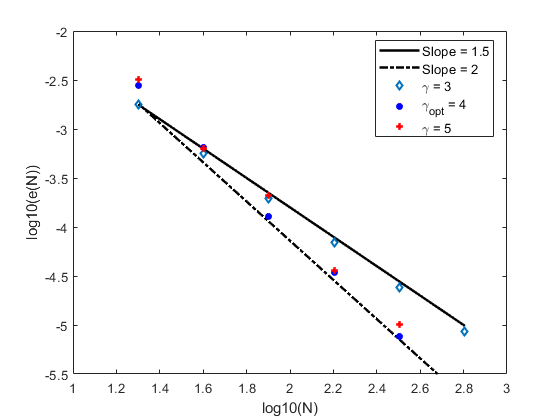} 
	\caption{ Temporal accuracy of \eqref{scheme: implicit Alikhanov TFAC} for $\alpha=0.8,\sigma=0.4$ and $\alpha=0.5,\sigma=0.5$.}\label{exam:numerical accuracy}
\end{figure}

We test the time accuracy with the  $L^2$ norm errors
$e(N):=\max_{1\le n\le N}\mynormb{\Phi^n-\phi^n}$.
By setting different grading parameters $\gamma>1$, the numerical results in Figure \ref{exam:numerical accuracy} are computed for $\sigma=0.4,\,\alpha=0.8$  and $\sigma=\alpha=0.5$, respectively.
As observed, the variable-step L2-1$_\sigma$ scheme \eqref{scheme: implicit Alikhanov TFAC} is of order $O(\tau^{\gamma \sigma})$ if the grading parameter $\gamma<\gamma_{\text{opt}}=2/\alpha$, while the optimal time accuracy $ O(\tau^{2})$ can be attained when the grading parameter $\gamma\ge\gamma_{\text{opt}}$.

\subsection{Simulation of coarsening dynamics}
\begin{example}
The coarsening dynamics of the TFAC model is examined with  
$\epsilon=0.05$. The initial condition is generated as $\Phi_{0}(\mathbf{x})=\mathrm{rand}(\mathbf{x}) $, where $\mathrm{rand}(\mathbf{x}) $ is uniformly distributed random number varying from $ -0.001 $ to 0.001 at each grid point.
\end{example}

\begin{figure}[htb!]
	\centering
	\includegraphics[width=2.0in]{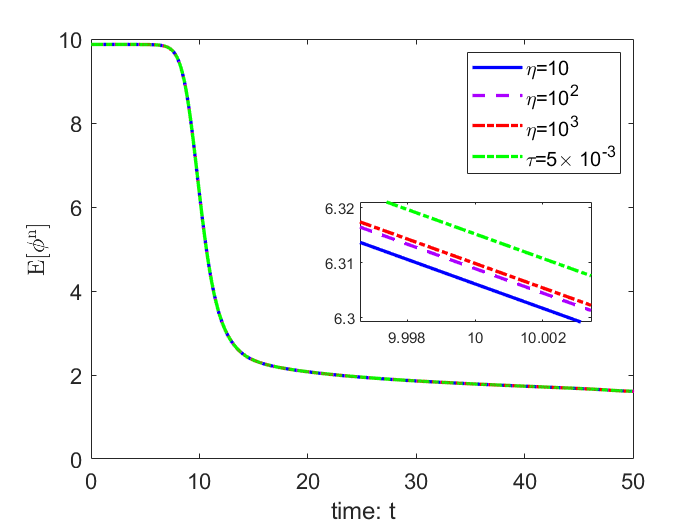}
	\includegraphics[width=2.0in]{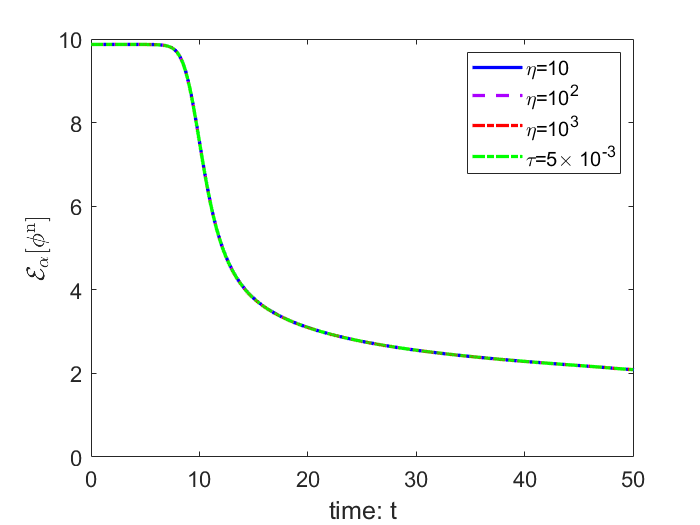}
	\includegraphics[width=2.0in]{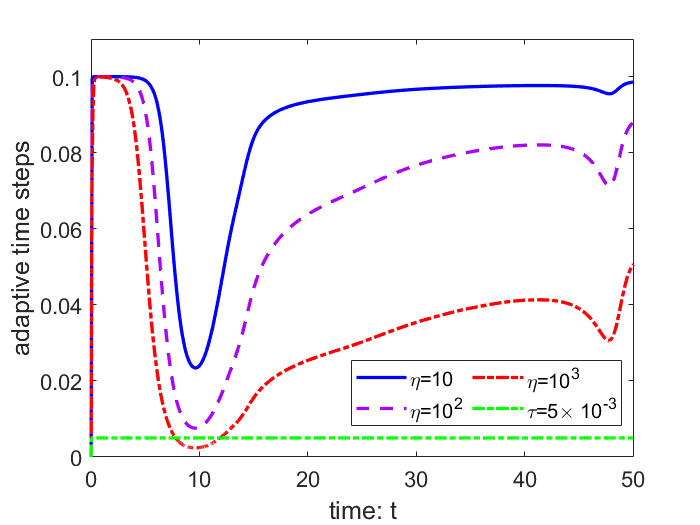}
	\caption{Energy curves from adaptive strategy for different parameters $\eta$.}
	\label{exam:TFAC compare eta}
\end{figure}

We use a $128^2$ grid to discretize the domain $(0,2\pi)^2$. The graded mesh $t_k=T_0(k/N_0)^{\gamma}$ with the settings $\gamma=3,N_0=30$ and $T_0=0.01$ are always applied in the initial part $[0,T_0]$, cf. \cite{JiLiao2020SimpleMPTFAC}.
In order to improve the computational efficiency of long-time numerical simulation and accurately capture the rapid changes of energy and numerical solution, we adopt the following adaptive time-stepping strategy to adjust the step size by the change rate of the solution \cite{Zhang2012AnAdaptiveCH,Huang2020ParallelACCH}
\begin{align*}
\tau_{ada}
=\max\big\{\tau_{\min},
\tau_{\max}/\varPi(\phi,\eta)\big\}
\quad \text{so that} \quad \tau_{n+1}=\max\{\tau_{ada}, r_{\star}\tau_n\},
\end{align*}
where $\tau_{\max}$ and $\tau_{\min}$ are the predetermined maximum and minimum size of time-steps, and 
$\varPi(\phi,\eta):=\sqrt{1+\eta\mynormb{\partial_\tau \phi^{n-\frac12}}^2}$ with a user parameter $\eta$ to be determined. 
We take $\tau_{\max}=10^{-1}$, $\tau_{\min}=10^{-3}$, $T_0=0.01$, $\gamma=3$ and consider three different parameters $\eta = 10 , 10^2 $ and $ 10^3$ to determine a suitable parameter $\eta$ where the fractional index $\alpha=0.7$, T=50 and the reference solution is computed by using the uniform time step $\tau=5\times 10^{-3}$. As seen in Figure \ref{exam:TFAC compare eta}, the value of parameter $\eta$ evidently influences on the adaptive sizes of time steps. The time-steps have the smallest fluctuation if $\eta=10^3$ and the energy curve is closer to the reference curve.

\begin{figure}[htb!]
\centering
\includegraphics[width=1.42in]{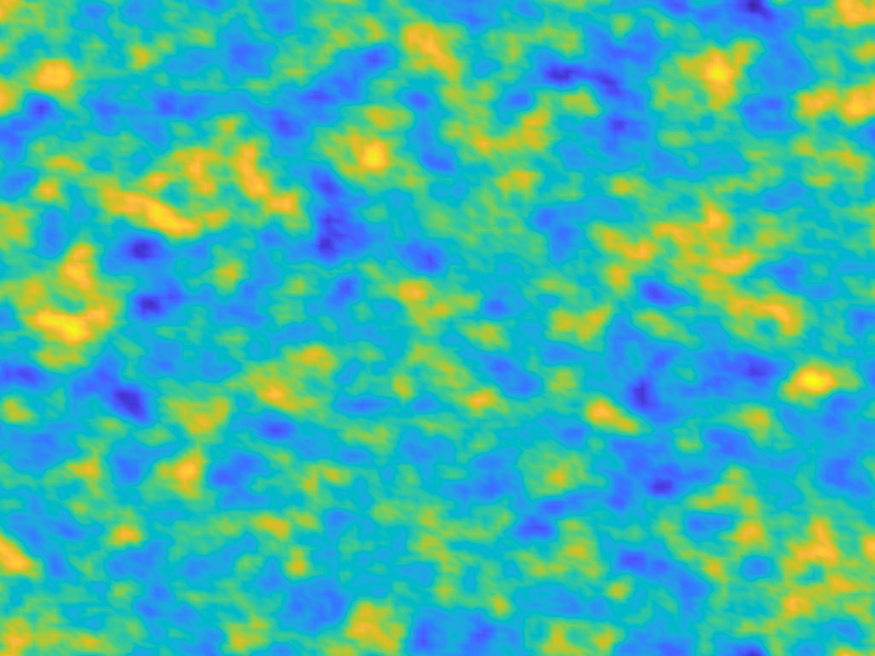}
\includegraphics[width=1.42in]{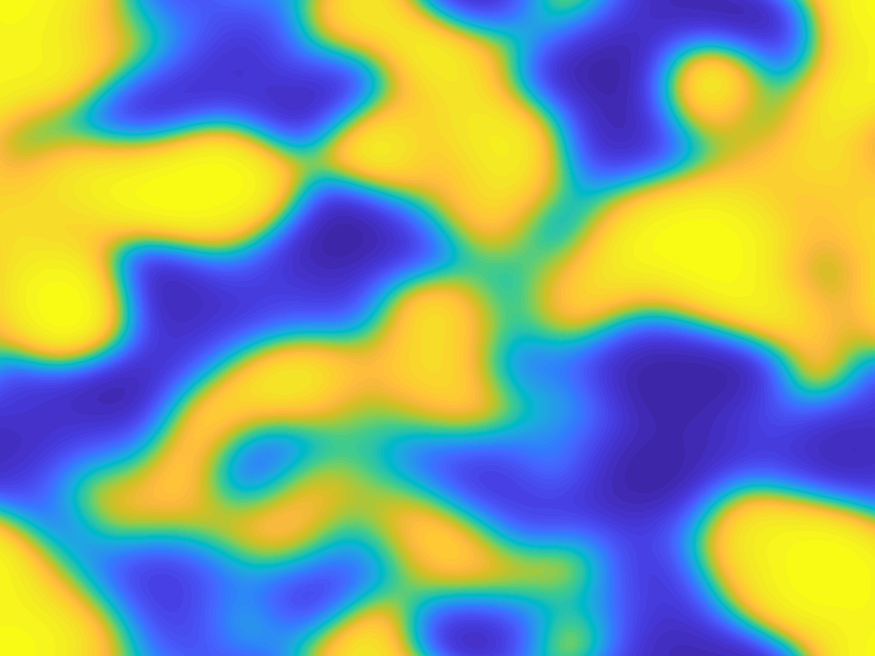}
\includegraphics[width=1.42in]{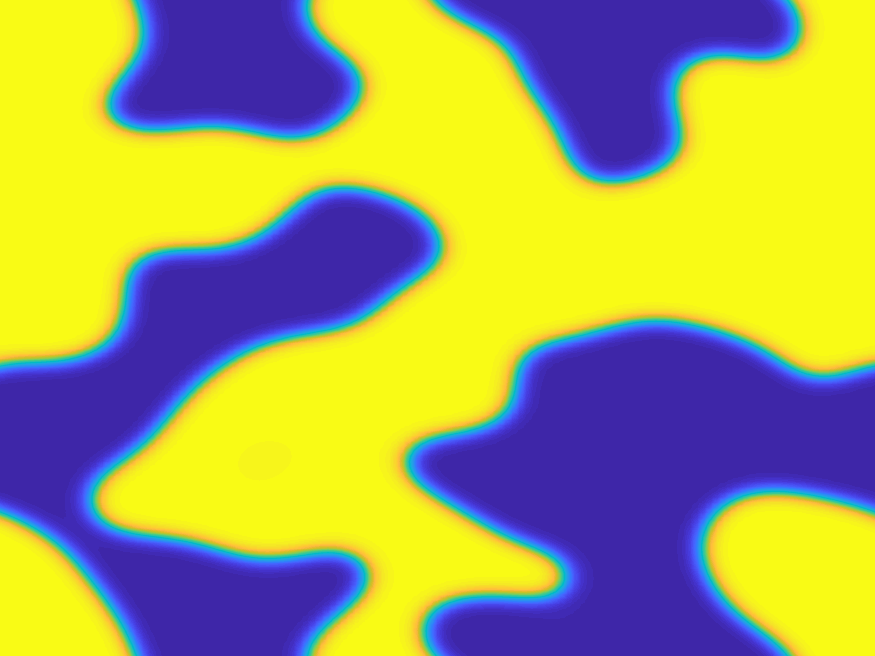}
\includegraphics[width=1.42in]{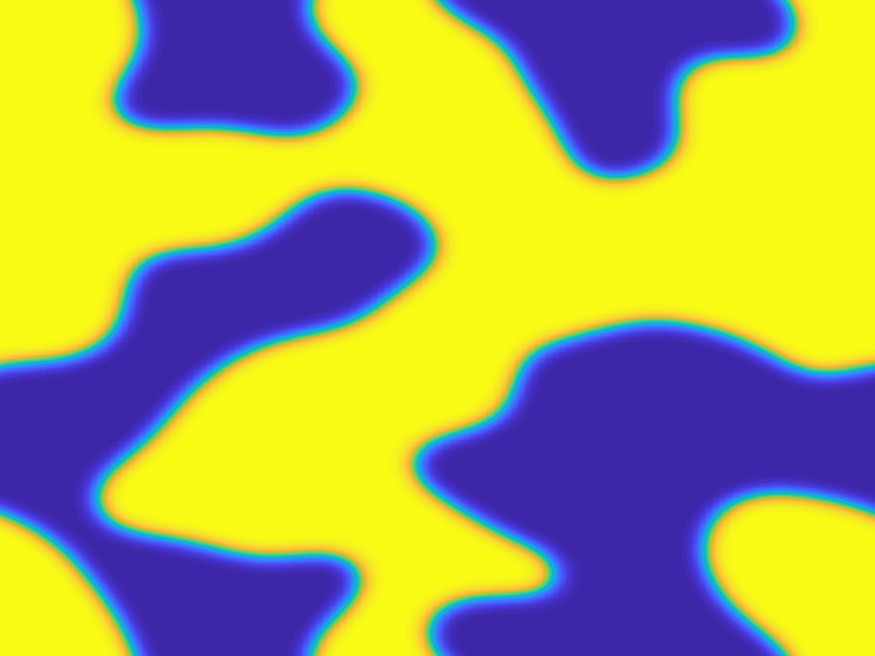}\\
\includegraphics[width=1.42in]{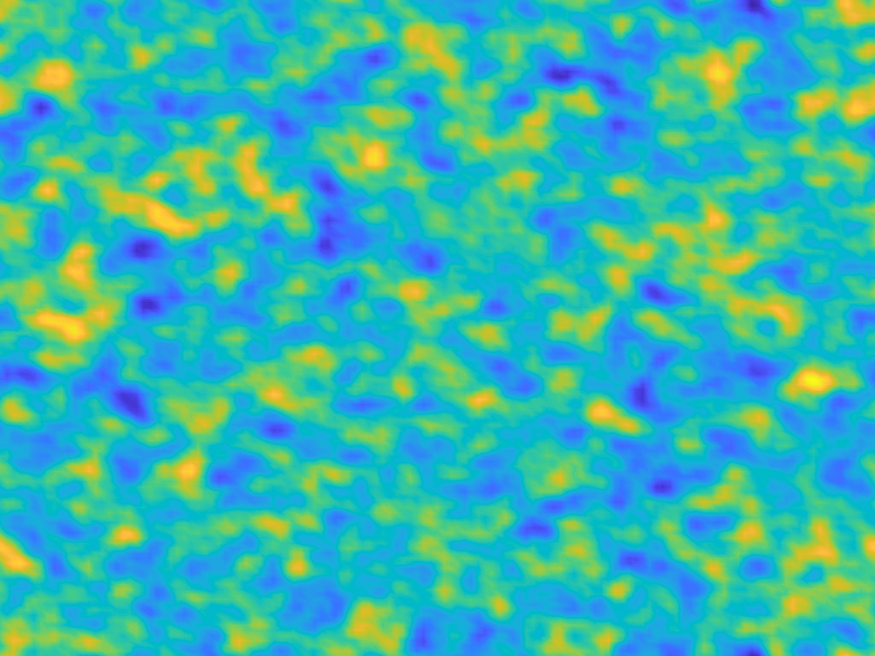}
\includegraphics[width=1.42in]{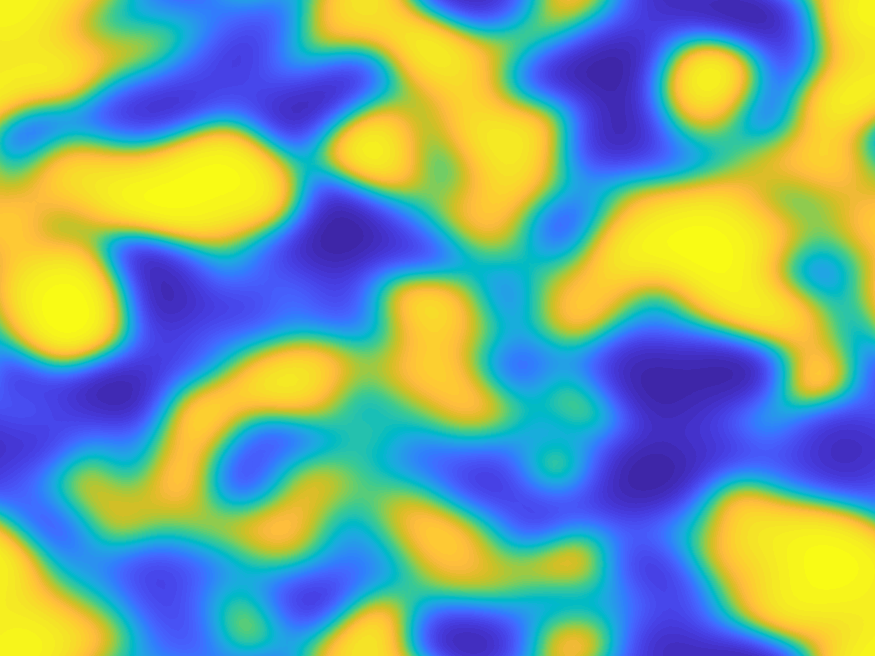}
\includegraphics[width=1.42in]{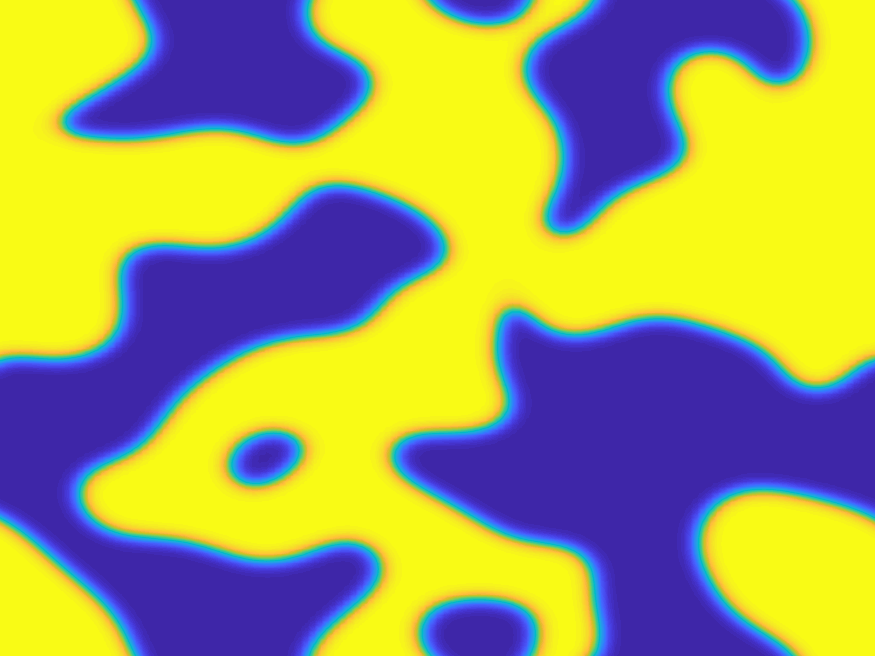}
\includegraphics[width=1.42in]{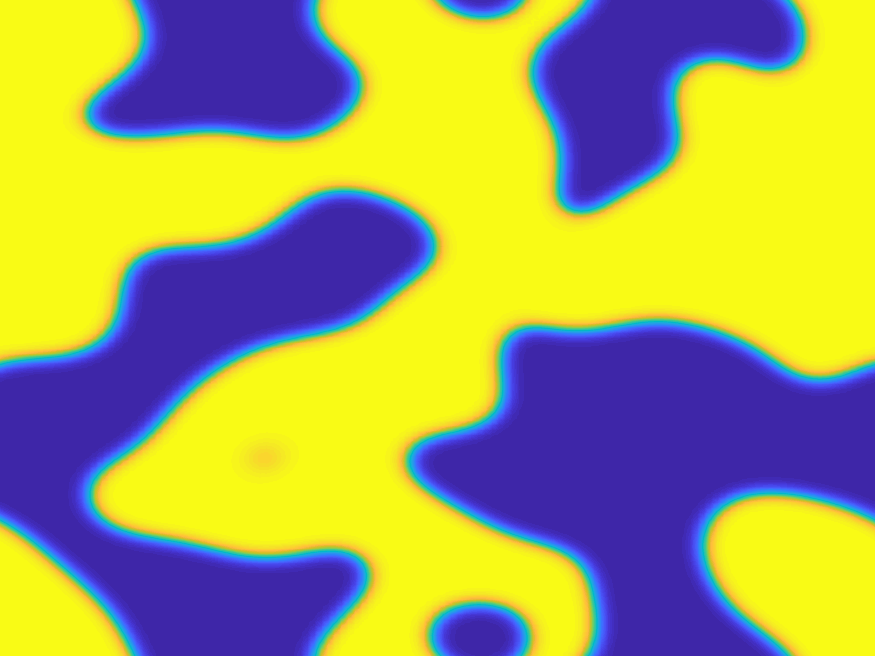}\\
\includegraphics[width=1.42in]{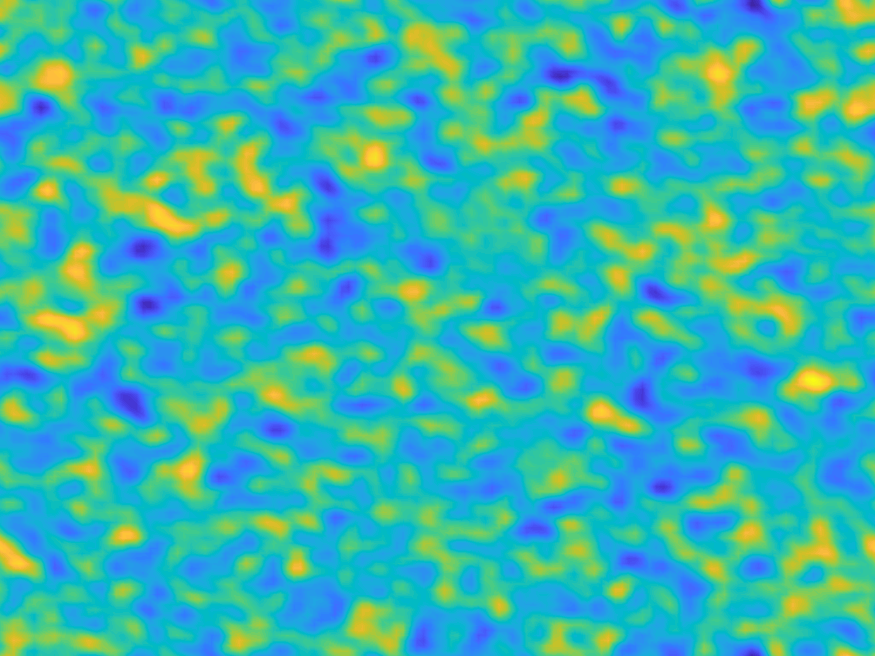}
\includegraphics[width=1.42in]{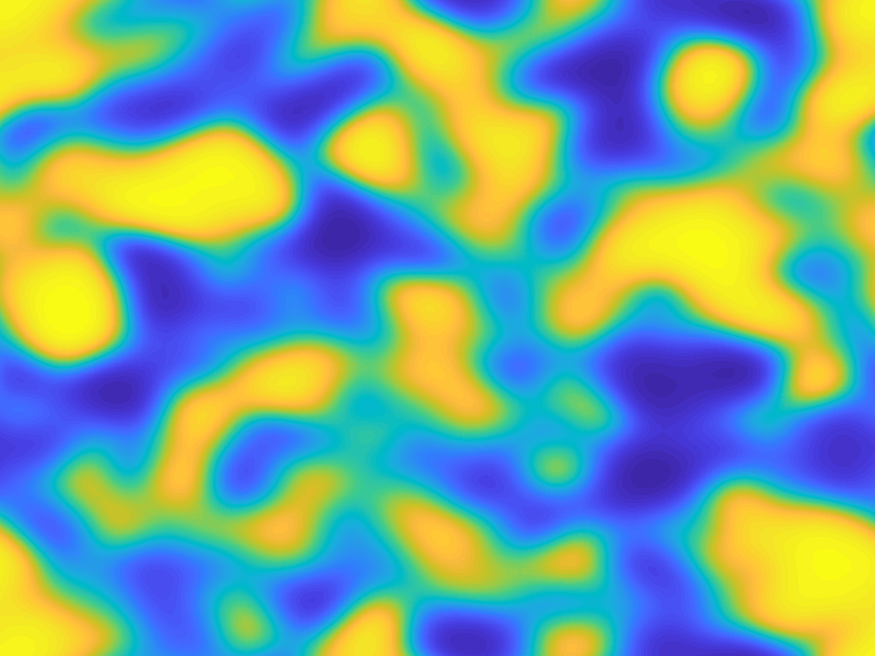}
\includegraphics[width=1.42in]{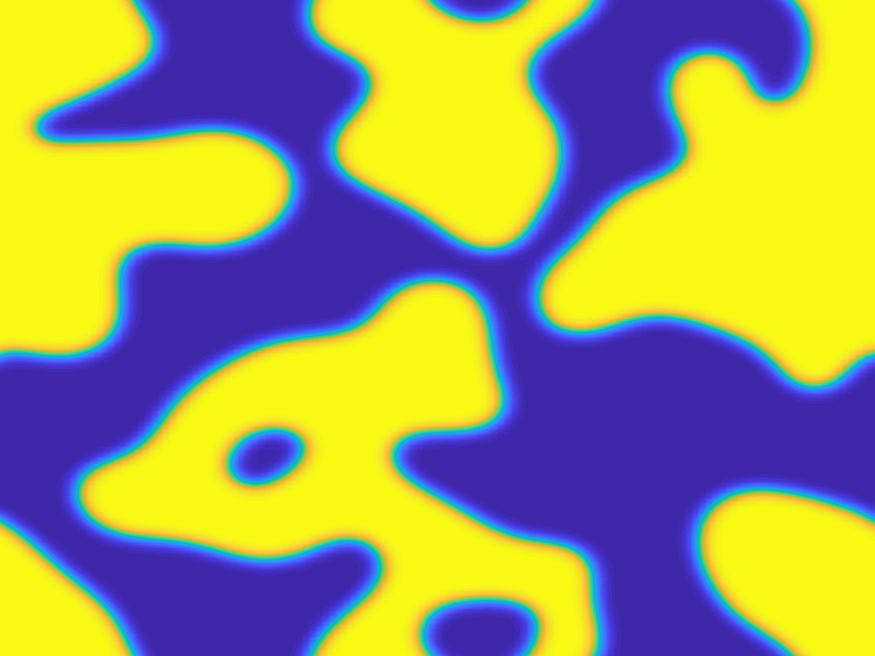}
\includegraphics[width=1.42in]{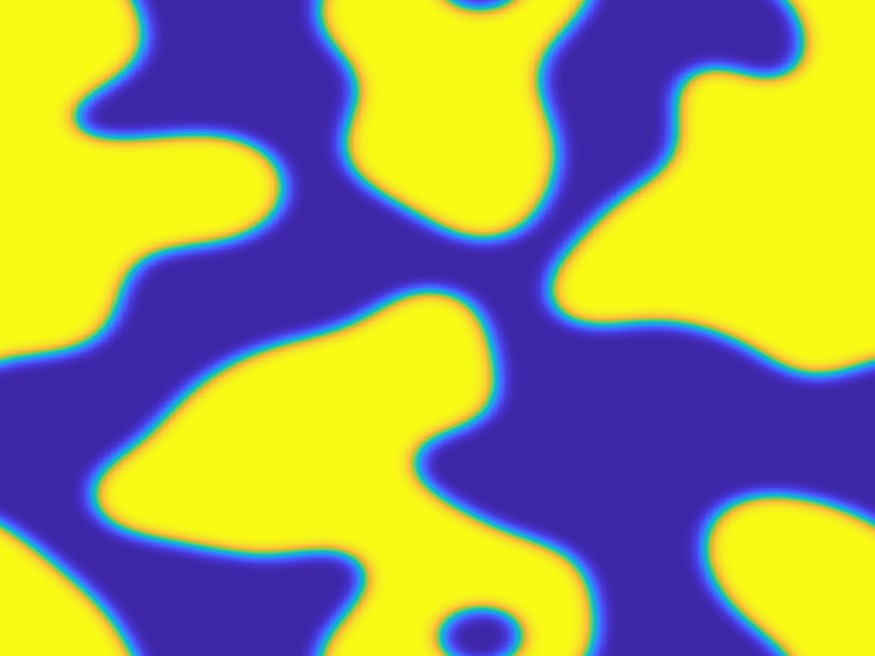}
\caption{Profiles of coarsening dynamics at $t=1, 10, 30, 50$ (from left to right) 
	for three different fractional orders $\alpha=0.4, 0.7, 0.9$ (from top to bottom).}
  \label{exam:TFAC dynamic snap}
\end{figure}

\begin{figure}[htb!]
	\centering
	\includegraphics[width=2.0in]{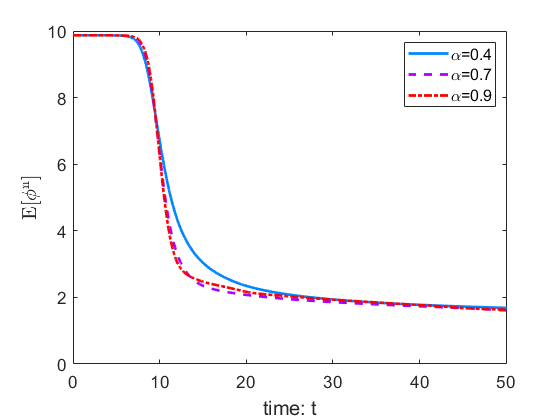}
	\includegraphics[width=2.0in]{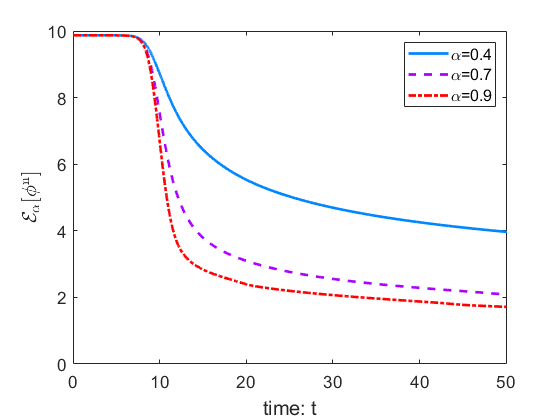}
	\includegraphics[width=2.0in]{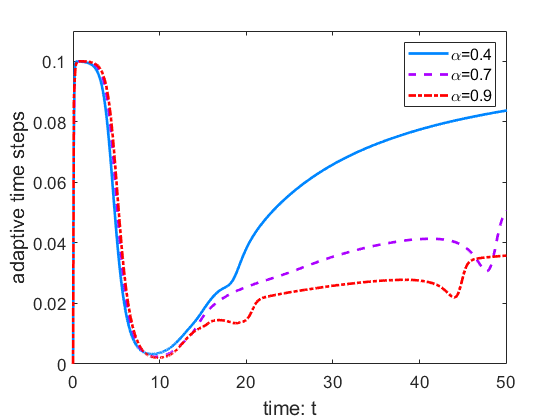}
	\caption{$E(t)$, $\mathcal{E}_\alpha(t)$ and adaptive time steps for the fractional orders $\alpha=0.4, 0.7, 0.9$.}
	\label{exam:TFAC dynamic energy}
\end{figure}

Next, by taking $\tau_{\max}=10^{-1}$, $\tau_{\min}=10^{-3}$ and the parameter $\eta=10^3$ in the adaptive time-stepping strategy, the  profiles of coarsening dynamics with different fractional orders $\alpha=0.4,0.7$ and 0.9 for the TFAC model are shown in Figure \ref{exam:TFAC dynamic snap}.
Snapshots are taken at time $t=1, 10, 30$ and 50, respectively.
We observe that the coarsening rates of TFAC model are dependent on the fractional order and the time period.
The calculated energies using the adaptive time steps for
the coarsening dynamics are depicted in Figure \ref{exam:TFAC dynamic energy}, which shows that the original energy and the modified energy for TFAC model with  $\alpha=0.4,0.7$ and 0.9 are both decreasing respect to time. And the  adaptive time steps strategy is successful in adjusting the time steps.

\begin{figure}[htb!]
	\centering
	\includegraphics[width=2.0in]{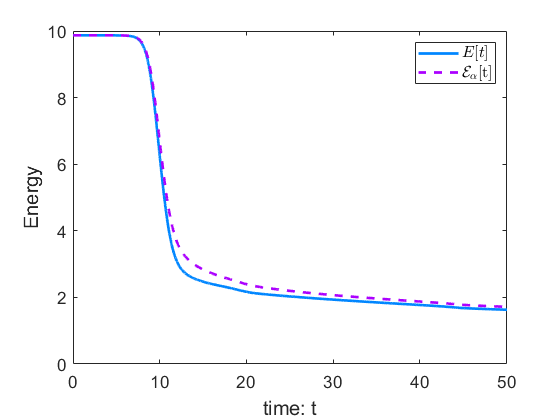}
	\includegraphics[width=2.0in]{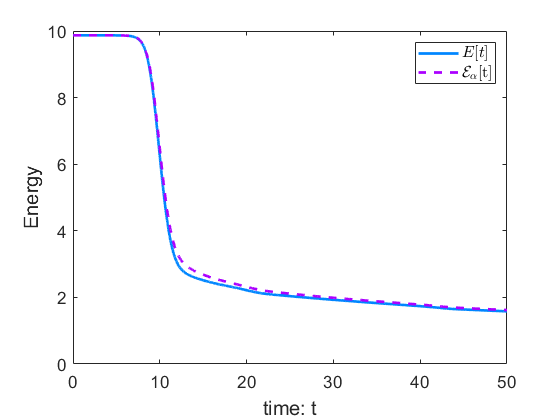}
	\includegraphics[width=2.0in]{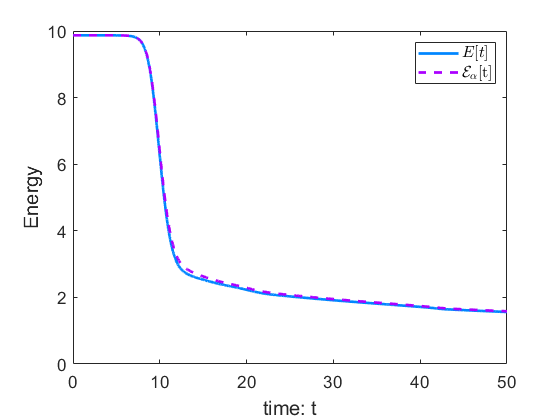}
	\caption{Comparison of $E(t)$ and $\mathcal{E}_\alpha(t)$
		for $\alpha=0.9, 0.95, 0.97$ (from left to right).}
	\label{exam:TFAC Asymptotic compatible energy}
\end{figure}

Figure \ref{exam:TFAC Asymptotic compatible energy} compares the curves of $E(t)$ and $\mathcal{E}_\alpha(t)$ for different fractional orders $\alpha=0.9, 0.95, 0.97$ to illustrate the asymptotic compatibility of the new discrete energy. As expected, the closer the time fractional index is to 1, the closer the modified energy $\mathcal{E}_\alpha(t)$ is to the original energy $E(t)$. 

\vspace{0.5cm}

\begin{description}	
	\item[\textbf{Funding}]   This work is supported by NSF of China under Grant Number 12071216.
	\item[\textbf{Data Availibility Statement}]  All data generated or analysed during this study are included in this published	article.
	\item[\textbf{Declarations}]  The authors declare that they have no conflict of interest.
			\item[\textbf{Compliance with Ethical Standards}]  This article does not contain any studies involving animals performed by any of the authors.
\end{description} 



\begin{thebibliography}{10}
	
	\bibitem{Alikhanov2015NewScheme}
	A.~A. Alikhanov.
	\newblock A new difference scheme for the time fractional diffusion equation.
	\newblock {\em J. Comput. Phys.}, 280:424-438, 2015.
	
	\bibitem{AllenCahn1979AMicroscopicTheory}
	S.~Allen and J.~Cahn.
	\newblock A microscopic theory for antiphase boundary motion and its
	application to antiphase domain coarsening.
	\newblock {\em Acta Metall.}, 27:1085-1095, 1979.
	
	\bibitem{DuYangZhou202TimeFractionalAllenCahn}
	Q.~Du, J.~Yang and Z.~Zhou.
	\newblock Time-fractional {Allen-Cahn} equations: Analysis and numerical
	methods.
	\newblock {\em J. Sci. Comput.}, 85:11, 2020.
	
	\bibitem{Feng2003NumericalAnalysisAllenCahn}
	X.~Feng and A.~Prohl.
	\newblock Numerical analysis of the {Allen-Cahn} equation and approximation
	for mean curvature flows.
	\newblock {\em Numer. Math.}, 94:33-65, 2003.
	
	\bibitem{HouXu2022sisc}
	D.~Hou and C.~Xu.
	\newblock Highly efficient and energy dissipative schemes for time fractional Allen-Cahn equation.
	\newblock {\em SIAM J. Sci. Comput.}, 43(5):A3305-A3327,2021.
	
	
	
	\bibitem{HouZhuXu2021L2SAVtfac}
	D.~Hou, H. Zhu and C.~Xu.
	\newblock Highly efficient schemes for time-fractional
	Allen–Cahn equation using extended SAV approach.
	\newblock {\em Numer. Algorithms}, 88:1077-1108, 2021.
		
	\bibitem{HouXu2022secondorderSAV}
	D.~Hou and C.~Xu.
	\newblock A second order energy dissipative scheme for time fractional $L^2$
	gradient flows using {SAV} approach.
	\newblock {\em J. Sci. Comput.}, 90(1):1-22, 2022.
	
	\bibitem{HouTangYang2017maximumprincipleCNTAC}
	T.~Hou, T.~Tang and J.~Yang.
	\newblock Numerical analysis of fully discretized {Crank-Nicolson} scheme for
	fractional-in-space {Allen-Cahn} equations.
	\newblock {\em J. Sci. Comput.}, 72(3):1214-1231, 2017.
	
	\bibitem{Huang2020ParallelACCH}
	J.~Huang, C.~Yang and Y.~Wei.
	\newblock Parallel energy-stable solver for a coupled {Allen-Cahn} and
	{Cahn-Hilliard} system.
	\newblock {\em SIAM J. Sci. Comput.}, 42:C294-C312, 2020.
	
	\bibitem{JiLiao2020SimpleMPTFAC}
	B.~Ji, H.-L. Liao and L.~Zhang.
	\newblock Simple maximum-principle preserving time-stepping methods for
	time-fractional {Allen-Cahn} equation.
	\newblock {\em Adv. Comput. Math.}, 46:37, 2020.
	
	\bibitem{JiLiaoGongZhang2020SISC}
	B. Ji, H.-L. Liao, Y. Gong, and L. Zhang. Adaptive second-order Crank-Nicolson timestepping
	schemes for time fractional molecular beam epitaxial growth models. SIAM J.
	Sci. Comput., 42:B738–B760, 2020.
	
	\bibitem{JiZhuLiao2022TFCH}
	B.~Ji, X.~Zhu and H-L. Liao.
	\newblock Energy stability of variable-step {L1}-type schemes for
	time-fractional {Cahn-Hilliard} model.
	\newblock {\em arXiv:2201.00920}, 2022.
	
	\bibitem{Jiang2017FastEvaluationCaputo}
	S.~Jiang, J.~Zhang, Z.~Qian and Z.~Zhang.
	\newblock Fast evaluation of the {C}aputo fractional derivative and its
	applications to fractional diffusion equations.
	\newblock {\em Commun. Comput. Phys.}, 21:650-678, 2017.
	
	 \bibitem{Karaa:2021}
	{S. Karaa}, {Positivity of discrete time-fractional operators 
		with applications to phase-field equations},
	{\em SIAM J. Numer. Anal.}, 59:2040-2053, 2021.
	
	
	\bibitem{Kopteva:2020}
	{N. Kopteva}, {Error analysis for time-fractional semilinear parabolic 
		equations using upper and lower solutions},
	{\em SIAM J. Numer. Anal.}, 58:2212-2234, 2020.
	
	
	\bibitem{KoptevaMeng:2020}
	{N. Kopteva and X. Meng}, 
	{Error analysis for a fractional-derivative parabolic problem 
		on quasi-graded meshes using barrier functions},
	{\em SIAM J. Numer. Anal.}, 58:1217-1238, 2020.
	
	\bibitem{LiWangYang:2017}
	Z. Li, H. Wang and D. Yang, A space-time fractional phase-field model with tunable sharpness
	and decay behavior and its efficient numerical simulation, 
	{\em J. Comput. Phys.}, 347: 20-38, 2017.
	
	\bibitem{LiaoLiuZhao:2022FBDF2}
	H.-L. Liao, N.~Liu and X.~Zhao.
	\newblock Asymptotically compatible energy of variable-step fractional BDF2 formula 
	for time-fractional Cahn-Hilliard model.
	\newblock {\em arXiv: 2210.12514v1}, 2022. 
	
	
	\bibitem{Liao2018DiscreteGronwallInequality}
	H.-L. Liao, W.~McLean and J.~Zhang.
	\newblock A discrete {G}r\"{o}nwall inequality with application to numerical
	schemes for subdiffusion problems.
	\newblock {\em SIAM J. Numer. Anal.}, 57:218-237, 2019.
	
	\bibitem{Liao2021SecondorderReactionsudiffusion}
	H.-L. Liao, W.~McLean and J.~Zhang.
	\newblock A second-order scheme with nonuniform time steps for a linear
	reaction-subdiffusion problem.
	\newblock {\em Commun. Comput. Phys.}, 30(2):567-601, 2021.
	
	\bibitem{Liao2020ASecondorderJCP}
	H.-L. Liao, T.~Tang and T.~Zhou.
	\newblock A second-order and nonuniform time-stepping maximum-principle
	preserving scheme for time-fractional {Allen-Cahn} equations.
	\newblock {\em J. Comput. Phys.}, 414:109473, 2020.
	
	\bibitem{LiaoTangZhou2021EnergyStableTFAC}
	H.-L. Liao, T.~Tang and T.~Zhou.
	\newblock An energy stable and maximum bound preserving scheme with variable
	time steps for time fractional {Allen-Cahn} equation.
	\newblock {\em SIAM J. Sci. Comput.}, 43(5):A3503-A3526, 2021.
	
	\bibitem{LiaoZhuWang:2021TFAC}
	H.-L. Liao, X.~Zhu and J.~Wang.
	\newblock An adaptive {L1} time-stepping scheme preserving a compatible energy
	law for the time-fractional {Allen-Cahn} equation.
	\newblock {\em Numer. Math. Theor. Meth. Appl.}, 15(4):1128-1146, 2022.
	
	\bibitem{LiuChengWangZhao2018TimefractionalACCH}
	H.~Liu, A.~Cheng, H.~Wang and J.~Zhao.
	\newblock Time-fractional {Allen-Cahn} and {Cahn-Hilliard} phase-field models
	and their numerical investigation.
	\newblock {\em Comput. Math. Appl.}, 76:1876-1892, 2018.
	
	\bibitem{QuanTangWang2022DecreasingUpperBound}
	C.~Quan, T.~Tang, B.~Wang and J.~Yang.
	\newblock A decreasing upper bound of energy for time-fractional phase-field
	equations.
	\newblock {\em arXiv:2202.12192v1}, 2022.
	
	\bibitem{Quan2020HowToDefine}
	C.~Quan, T.~Tang and J.~Yang.
	\newblock How to define dissipation-preserving energy for time-fractional
	phase-field equations.
	\newblock {\em CSIAM-AM}, 1:478-490, 2020.
	
	\bibitem{QuanWu2022stabilityL21sigma}
	C.~Quan and X.~Wu.
	\newblock On stability and convergence of {L2-1$_\sigma$} method on general
	nonuniform meshes for subdiffusion equation.
	\newblock {\em arXiv:2208.01384v1}, 2022.
	
	\bibitem{QuanWu:2022H1stabilityL2}
	C. Quan and X. Wu,
	H1-stability of an L2 method on general nonuniform meshes for subdiffusion equation, 
	{\em arXiv:2205.06060v1}, 2022.
	
	\bibitem{Tang2018OnEnergyDissipation}
	T.~Tang, H.~Yu and T.~Zhou.
	\newblock On energy dissipation theory and numerical stability for
	time-fractional phase field equations.
	\newblock {\em SIAM J. Sci. Comput.}, 41:A3757-A3778, 2019.
	
	\bibitem{Zhang2012AnAdaptiveCH}
	Z.~Zhang and Z.~Qiao.
	\newblock An adaptive time-stepping strategy for the {Cahn-Hilliard} equation.
	\newblock {\em Commun. Comput. Phys.}, 11:1261-1278, 2012.
	
	\bibitem{ZhaoChenWang:2019}
	J. Zhao, L. Chen and H. Wang, On power law scaling dynamics for time-fractional phase field
	models during coarsening. {\em Commun. Nonlinear Sci.}, 70:257-270, 2019.
	
\end{thebibliography}

\end{document}